# Hamiltonian Graphs and the Traveling Salesman Problem

Dhananjay P. Mehendale
Sir Parashurambhau College, Tilak Road, Pune 411030,
India

## Abstract

A new characterization of Hamiltonian graphs using *f*-cutset matrix is proposed. A new exact polynomial time algorithm for the traveling salesman problem (TSP) based on this new characterization is developed. We then define so called ordered weighted adjacency list for given weighted complete graph and proceed to the main result of the paper, namely, the exact algorithm based on utilization of ordered weighted adjacency list and the simple properties that any path or circuit must satisfy. This algorithm performs checking of sub-lists, containing (*p*-1) entries (edge pairs) for paths and *p* entries (edge pairs) for circuits, chosen from ordered adjacency list in a well defined sequence to determine exactly the shortest Hamiltonian path and shortest Hamiltonian circuit in a weighted complete graph of *p* vertices. The procedure has intrinsic advantage of landing on the desired solution in quickest possible time and even in worst case in polynomial time. A new characterization of shortest Hamiltonian tour for a weighted complete graph satisfying triangle inequality (i.e. for tours passing through every city on a realistic map of cities where cities can be taken as points on a Euclidean plane) is also proposed. Finally, we propose a classical algorithm for unstructured search, three new quantum algorithms for unstructured search, which exponentially speed up the searching ability in the unstructured database and one quantum algorithm for solving a K-SAT problem and indicate its effect on traveling salesman problem and other NP-Complete problems.

**1. Characterization of Hamiltonian Graphs using *f*-cutset Matrix:** Let *G* be a (p, q) graph, i.e. a graph on p points (vertices) and q lines (edges) with the following vertex set $V(G)$ and edge set $E(G)$ respectively:

$$V(G) = \{ v_1, v_2, \cdots, v_p \} \text{ and}$$

$$E(G) = \{ e_1, e_2, \cdots, e_q \}$$

Let $A_G = [a_{ij}]_{p \times p}$ denotes the adjacency matrix of *G*.

By choosing a spanning tree in the given connected graph one can construct the fundamental cutest matrix associated with this choice of tree in the form

$$C_f = [C_c : I_{p-1}]$$

where the last $(p-1)$ columns forming the identity matrix correspond to $(p-1)$ branches of the spanning tree, and the first $(q-p+1)$ columns forming matrix $C_c$ correspond to the chords. The presence of entry "1" in a column and row of $C_f$ depicts the presence of the edge represented by that column in the *f*-cutest represented by that row. Each row in $C_f$ is a fundamental cutest vector. The rank of $C_f$ is $(p-1)$ and these fundamental cutest vectors form the vector space basis of the cutest space, which is subspace of the vector space associated with the graph. If we develop an algorithm which selects edges from the fundamental cutsets (*f*-cutsets) to form a Hamiltonian circuit then it is clear to see that

1) Since every vertex must get incorporated in every Hamiltonian circuit and only once, so, every f-cutset must contribute positive and even number of edges.
2) Since $C_f$ partitions into $C_c$ and Identity matrix, $I_{p-1}$, and since every row of an identity matrix contains at most one nonzero entry, so every *f*-cutset must contribute at least one chord to every Hamiltonian circuit in order to maintain evenness and positivity of the number of edges chosen on the corresponding *f*-cutset. Thus, the presence of a branch from every *f*-cutset in a Hamiltonian circuit is not necessary but the presence of at least one chord from every *f*-cutset is a must for a Hamiltonian circuit formed by chords and branches.

**Theorem 1.1:** Let *G* be a (p, q) graph. *G* is Hamiltonian if and only if we can select even number (at least two) of edges on each row representing an $f$–cutset to form a connected graph such that at least one of the selected edge is a chord and the total count of thus selected edges (chords + branches) is *p*.

**Proof:** Let G be a Hamiltonian graph. So, there is a Hamiltonian circuit in G. So, there exists a tree which is Hamiltonian path. Take this tree which is a Hamiltonian path containing $(p-1)$ edges. Form *f*-cutset matrix for this tree. Among the chords represented by columns of $C_c$ there will exist a chord (which is actually the remaining part of the Hamiltonian circuit to be added to the tree equal to the Hamiltonian path to complete the Hamiltonian circuit) such that its corresponding column will be entirely made up of units. Thus, we take this chord and the corresponding branch on every *f*-cutset, in effect, even (= 2) number of edges are selected from each row representing an $f$ –cutset and the total count of thus selected edges (1 chord + $(p-1)$ branches) is equal to *p*, as desired forming a connected graph.

Suppose we have formed a subgraph of the given graph containing at least two edges on each row corresponding to every $f$ –cutset to form a connected subgraph such that at least one is a chord, so that, in effect even (>0) number of edges get selected from each row representing an $f$ –cutset and the total count of thus selected edges (chords + branches) is *p*. We show that this subgraph must be a Hamiltonian circuit. It can be easily seen that such a graph can have only two possibilities: Either it is a circuit of length *p*, or a subcircuit of length smaller than *p* with some incident (one or more) paths to the vertices of this subcircuit. In the first case nothing to prove. In the other case some *f*-cutset among the *f*-cutsets to which the edge incident on the pendant point of the path belongs, may be as a branch or a chord, must be contributing odd number of edges since there is no provision to reach and go away in the subgraph from this pendant point, a contradiction to the data.

□

**Some Interesting Observations:** (1) Search whether there is a column vector in $C_c$ entirely containing units (one unit in each row) i.e. in total $(p-1)$ units. In this case, this edge (chord) represented by the column of units and the edges representing all branches together sum up to in all $1+(p-1)=p$, edges forming a Hamiltonian circuit. Thus, the graph will be Hamiltonian.

(2)Search whether there are some two columns vectors in $C_c$ such that there is a unit in some rows in that column corresponding to a chord and there is a unit in the remaining rows in the column corresponding to the other chord among the total $(p-1)$ rows with exactly one overlap, i.e. there exists exactly one row which contains units in the chosen columns corresponding to both the chords. In this case, the two chosen edges (chords) represented by the two columns (determined as above with exactly one overlap) of $C_c$ and the edges representing all branches except the branch defining an *f*-cutset for which the corresponding row contains units in both the columns corresponding to the two chosen chords, together sum up to in all $2+(p-2)=p$, edges forming a Hamiltonian circuit. Thus, the graph will be Hamiltonian.

(3)Continuing on these lines, search if there exist some $k$ columns vectors in $C_c$ in which the units are distributed among these determined columns of $C_c$ such that there will exist in all $(k-1)$ overlaps as mentioned in (2), such that the $k$ chosen edges (chords) represented by the $k$ determined columns of $C_c$ and the edges representing all branches except the branches defining an *f*-cutset for which the corresponding row contains units in even number of columns corresponding to the $k$ chosen chords( i.e. there is overlap), and together sum up to in all $k+(p-k)=p$, edges forming a Hamiltonian circuit. Thus, the graph will be Hamiltonian.

□

We now state few definitions and develop a polynomial time algorithm using them:

**Definition 1.1:** A **lattice** is a rectangular array of dots made up of some *m* rows and *n* columns.

**Definition 1.2:** A **lattice-cutset-graph** associated with a $f$–cutset matrix is a graph obtained from the lattice of the size of $f$–cutset matrix obtained by treating those dots in the array as vertices of this new graph where, at the same place entry "1" is present in the $f$–cutset matrix, then joining these vertices in each column by vertical edges so that a path is formed along each column and further joining by horizontally going edges connecting all the vertices in a row to each other.

**Definition 1.3:** An **induced-cutset--tree** is a tree in the usual sense in the lattice-cutset-graph (and **not in the original graph**) such that its vertices belong to the first $(q-p+1)$ columns of lattice-cutset-graph and there exists at least one vertex belonging to every row of vertices of the lattice-cutset-graph and when some columns of vertices are chosen we should take all the vertical edges belonging to those columns and all the horizontally going edges joining the vertices in succession belonging to same row in the chosen columns.

**Definition 1.4:** An **extended-induced-cutset-tree** is the one obtained by induced-cutset-tree by adding horizontally going edges such that one horizontal edge is to be added connecting last (in fact, any vertex will do) vertex present in that row (belonging to the $C_c$ part) to the vertex belonging to the last $(p-1)$ columns ( belonging to the Identity part) in the same row only when the number of vertices that get incorporated in that row for the induced-cutset-tree are odd in number (due to the odd number of edges got selected in the $C_c$ part with respect to that row).

A simpler algorithm will be given below which just consists of finding an induced-cutset-tree and then an extended-induced-cutset-tree from this induced-cutset-tree.

**Remark 1.1:** Forming an induced-cutset-tree and then the extended-induced-cutset-tree automatically takes care of the important requirements in the above given observations of achieving exactly $(k-1)$ overlaps when $k$ chords are chosen, so that the important equation to be satisfied by the chosen chords and branches to form a Hamiltonian circuit, namely, chords + branches = $p$ is satisfied automatically.

**Steps to construct a lattice-cutset-graph and an extended-induced-cutset-tree as its subgraph from *f* –cutset matrix:**

(1) Form *f* –cutset matrix, $C_f = [C_c : I_{p-1}]$.
(2) Form a lattice of size (rows and columns of dots) equal to the size of the *f* –cutset matrix, made up of dots and then proceed to form lattice-cutset-graph by taking the dots with entry “1” in the same place in the associated fundamental cutest matrix and then by connecting the appropriate edges as per the above definition.
(3) Form an induced-cutset-tree from edges in the first $(q-p+1)$ columns (forming matrix $C_c$ of the *f* –cutset matrix that correspond to the chords) of the lattice-cutset-graph by choosing appropriate columns such that at least one vertex from each row of vertices gets incorporated.
(4) Count the number of vertices belonging to each row of this tree and when the number of vertices contained in the row (corresponding to the $C_c$ part of the *f* –cutset matrix) are odd in number then extend the tree by joining last (right most) vertex in the tree to the vertex in the same row and belonging to the last $(p-1)$ columns (corresponding to the identity matrix representing the $(p-1)$ branches of the $f$ –cutset matrix) by a new edge to obtain an extended-induced-cutset-tree.

**Example 1.1:** Consider the following graph, *H* say.

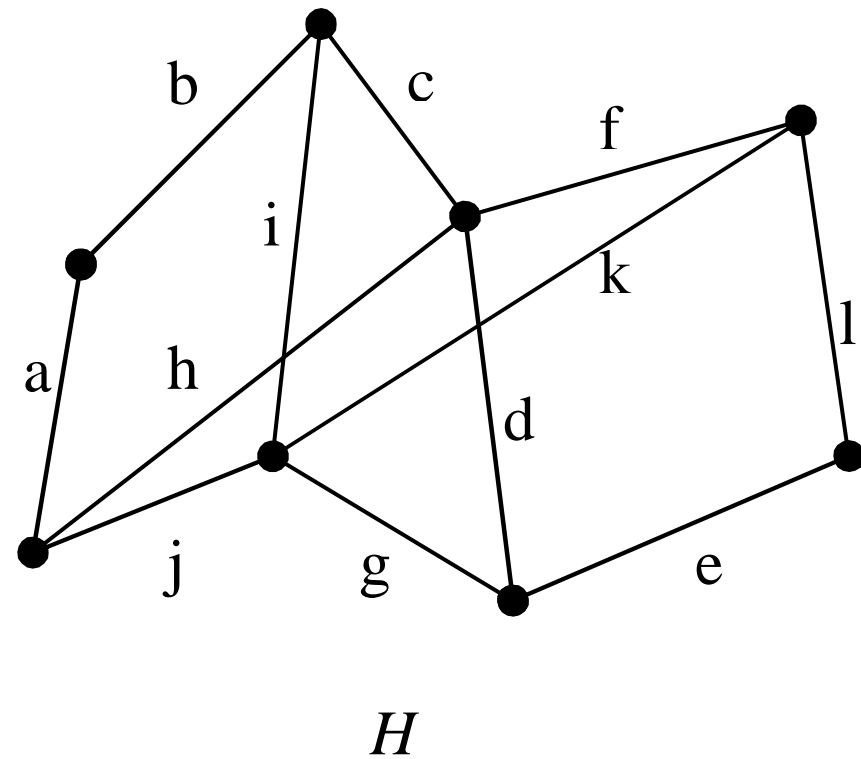

*H*

We take the spanning tree, $T$ , formed by edges {a, b, c, d, e, f, g} of this graph.

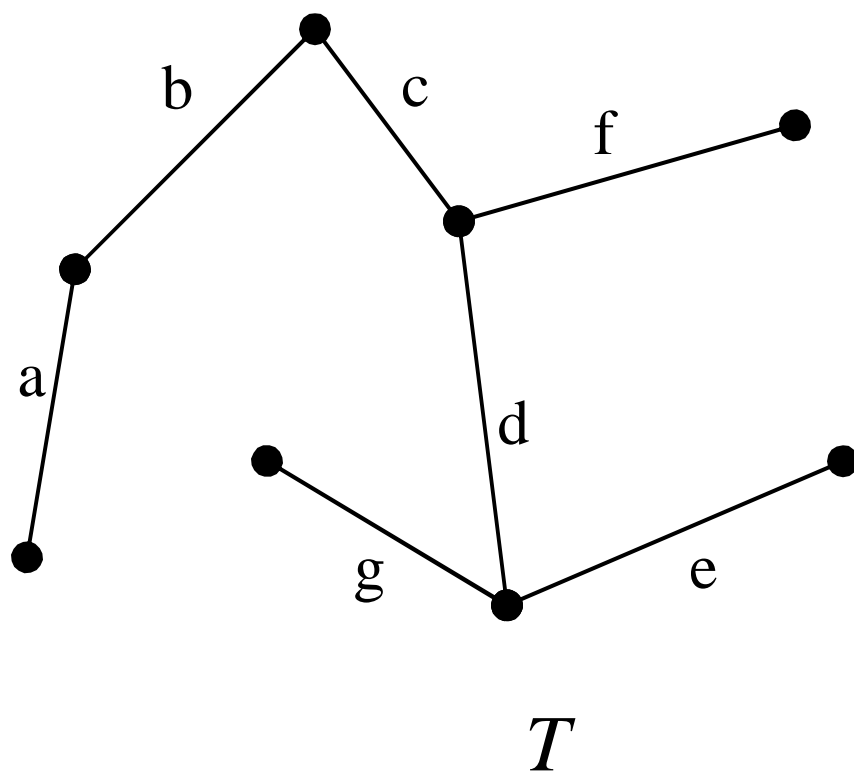

$T$

Then the $f$ –cutset matrix, $C_f$ can be expressed as follows, where the first column represents the labels of the $f$ –cutsets while the first row represents the labels of chords and branches. The $f$-cutset matrix, $C_f = [C_c : I_{p-1}]$
for the present case can be written as

$$
C_f = \begin{bmatrix} & h & i & j & k & l & a & b & c & d & e & f & g \\ c_1 & 1 & 0 & 1 & 0 & 0 & 1 & 0 & 0 & 0 & 0 & 0 & 0 \\ c_2 & 1 & 0 & 1 & 0 & 0 & 0 & 1 & 0 & 0 & 0 & 0 & 0 \\ c_3 & 1 & 1 & 1 & 0 & 0 & 0 & 0 & 1 & 0 & 0 & 0 & 0 \\ c_4 & 0 & 1 & 1 & 0 & 1 & 0 & 0 & 0 & 1 & 0 & 0 & 0 \\ c_5 & 0 & 0 & 0 & 0 & 1 & 0 & 0 & 0 & 0 & 1 & 0 & 0 \\ c_6 & 0 & 0 & 0 & 1 & 1 & 0 & 0 & 0 & 0 & 0 & 1 & 0 \\ c_7 & 0 & 1 & 1 & 0 & 0 & 0 & 0 & 0 & 0 & 0 & 0 & 1 \end{bmatrix}
$$

Consider the induced-cutset-trees, $T_1$ and $T_2$, containing chords {*j*, *l*} and chords {*h*, *i*, *l*} respectively, as shown below. An extended-induced-cutset-tree, $ET_1$ is also shown as an example. The first row of alphabets depicts the used **edge-labels**:

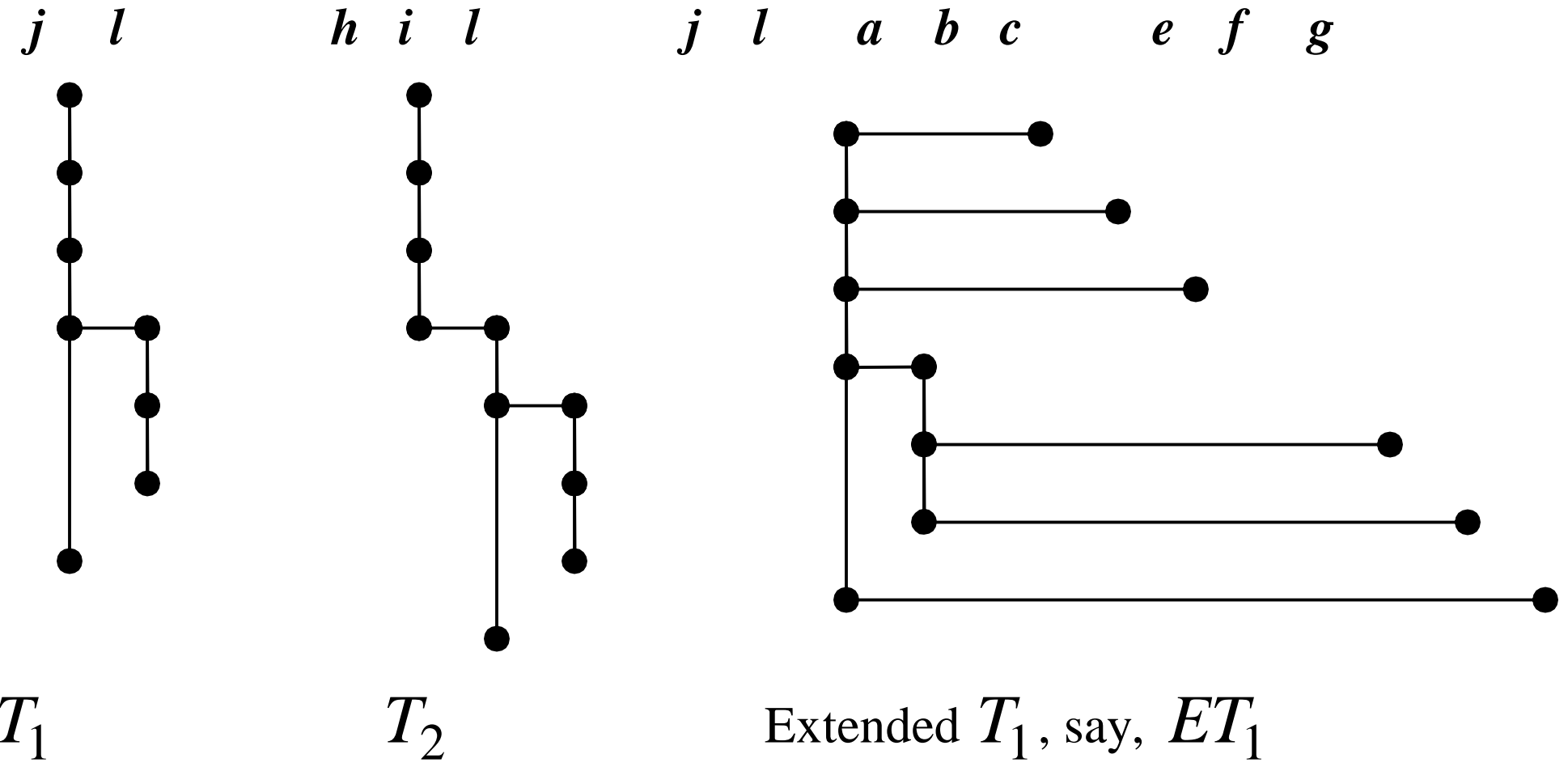

$T_1$ $T_2$ Extended $T_1$, say, $ET_1$

Now,

(i) Using chords {*j*, *l*} we form a tree, $T_1$, formed in $C_c$ and append the appropriate branches {*a*, *b*, c, *e*, *f* , *g*} leading to formation of extended tree, $ET_1$. This leads to the following Hamiltonian circuit:

$$a \rightarrow b \rightarrow c \rightarrow f \rightarrow l \rightarrow e \rightarrow g \rightarrow j \rightarrow a$$

(ii) Using the chords {*h*, *i* , *l*} we can see that we have a tree, $T_2$, containing formed in $C_c$ and appending the appropriate branches {*a*, *b*, *e*, *f* , *g*} we can get an extended tree, say $ET_2$, which will lead to the following Hamiltonian circuit :

$$a \rightarrow b \rightarrow i \rightarrow g \rightarrow e \rightarrow l \rightarrow f \rightarrow h \rightarrow a.$$

**Theorem 1.2:** An extended-induced-cutset-tree forms a Hamiltonian circuit if the total count of the chosen edges represented by columns of lattice-cutset-graph (chords + branches) is equal to *p*.

**Proof:** Straightforward. □

**Remark 1.2:** The edges $i$ and $j$ together cannot belong to any Hamiltonian circuit of graph *H* because we can't form an extended-induced-cutset-tree along with other edges such that the count of the selected edges (chords + branches) = *p*.

**2. The Traveling Salesman Problem (TSP):** This well-known problem asks for an efficient (polynomial time) algorithm to find shortest Hamiltonian circuit (or cycle), i.e. the one with smallest weight sum of its edges in a weighted complete graph. No efficient exact algorithm for this problem is known. It has been shown in the literature that the problem of finding Hamiltonian path between two pre-specified vertices, or that of finding Hamiltonian circuit, or that of finding shortest Hamiltonian tour for the traveling salesman etc. are all belong to the large class of **NP complete** problems (Page 234, Theorem 8.9 [2]). In fact it has been further shown that when the triangle inequality is not satisfied the traveling salesman problem is non-approximable unless NP complete problems have polynomial time solutions [3].

There exist some well-known efficient heuristic algorithms for TSP. The simplest one is the so called nearest-neighbor-method [4] with performance guarantee, $\alpha = \frac{1}{2}([\ln(n)]+1)$, and the efficient one is the so called minimum-weight- matching-algorithm [5] with performance guarantee $\alpha < \frac{3}{2}$.

In order to initiate the discussion we begin with a possible modification in the so called nearest-neighbor-method. The nearest-neighbor-method starts from a vertex $v_i$ in a weighted complete graph and select an edge among the edges the adjacent vertex $v_j$ such that the weight of the edge $e_k = (v_i, v_j)$ is minimum among the edges emerging from the vertex $v_i$. It continues with the same criterion for the incorporation of the next edge from vertex $v_j$ to a new vertex (i.e. not already visited one till all the vertices are exhausted and one has to select now the only left out choice $v_k$ to $v_i$ to complete the formation of the Hamiltonian circuit). Thus, to new nearest neighbor of $v_j$. Since the decision for choosing edge in this algorithm is based on purely **local** considerations, i.e. the selection is made which is locally best; there is no guarantee of this algorithm of attaining a good Hamiltonian circuit.

We propose below a modification in this algorithm which will make it **somewhat global** and thus will improve the chance of getting better performance. Let us take the given weighted complete graph as a symmetric digraph and each edge as two directed edges of same weight directed in opposite directions. It is clear to see that when one selects an edge, say $(v_i, v_j)$, then one cannot select any other edge emerging from $v_i$ and any edge entering in $v_j$. Thus, when one selects edge $(v_i, v_j)$ the weight sum of (other) edges emerging from $v_i$ and weight sum of (other) edges entering in $v_j$, other than edge $(v_i, v_j)$ gets excluded. Let us denote the weight-inclusion-adjacency- matrix by $WIA = [w^i_{jk}]_{p \times p}$, where $w^i_{jk}$ = weight of the edge $(v_j, v_k)$ that get included in the weight sum of edges of a Hamiltonian circuit when one selects the edge $(v_j, v_k)$ while forming that Hamiltonian circuit. Similarly, let us define the weight-exclusion- adjacency- matrix, $WEA = [w^e_{jk}]_{p \times p}$, where $w^e_{jk}$ = weight sum of the edges that get excluded while one selects the edge $(v_j, v_k)$ as an edge for a Hamiltonian circuit. Thus, $w^e_{jk} = \sum_{l \neq j} w^i_{lk} + \sum_{m \neq k} w^i_{jm}$

(As a further modification, we can add in $w^e_{jk}$ the weight sum of those edges which form a subcircuit with earlier selected edges and so can't be part of a Hamiltonian circuit. To keep the things simple here we do not take into consideration this further modification.)

**2.1 Modified-Nearest-Neighbor-Method:** In the nearest-neighbor-method there is only one criterion that is followed: the minimum weight nearest neighbor is selected to join, in succession by starting from some vertex, till one forms the Hamiltonian circuit. To make this algorithm somewhat global one follows the same selection method by imposing two criteria:

**Algorithm 2.1.1:** (1) Using the given $WIA = [w^i_{jk}]_{p \times p}$, construct the weight-exclusion-adjacency-matrix, $WEA = [w^e_{jk}]_{p \times p}$.

(2) Select an edge which obeys two criteria in any order:
  (i) The weight that gets included due to this selection is minimum.
  (ii) The weight that gets excluded due to this selection is maximum.

(3) Continue on these lines till one gets the desired Hamiltonian circuit.

□

Note that observing the criteria (i), (ii) in (2) **in any order and at any stage of the selection** it is possible to carry out the nearest- neighbor-algorithm, i.e. condition (i) is observed first and then (among the choices) condition (ii) is observed next, or, vice versa, and at any stage of the selection, i.e. we can change the order of conditions to be observed even at any intermediate stage of the algorithm.

**2.2 A Heuristic for TSP using Contractions:** The following heuristic is certainly an improvement over the usual nearest-neighbor-method because it eliminates the restriction of choosing only adjacent edges in succession imposed in the nearest-neighbor-method:
**Algorithm 2.2.1:** (1) Choose an edge among the edges with smallest weight in the given weighted complete graph on $p$ points and contract it. Keep the record of the contracted edge. This leads to formation of weighted complete graph on $p-1$ points.
(2) Repeat the procedure in step (1) for the resulted weighted complete graph on $p-1$ points till we reach to a simple graph with single vertex.
(3) Build the Hamiltonian circuit using the contracted edges whose record has been kept at every stage of contraction.

□

**Remark 2.1:** In the nearest-neighbor-method we select smallest edge among the adjacent edges but in the above given heuristic algorithm 2.4 each time we select the smallest available edge, not necessarily adjacent one, at each stage of selection. This idea of contraction can be taken up also in the modified-nearest-neighbor-method which will produce a complete digraph at each time, may be asymmetric in weight after some iterations, to continue the same procedure till a Hamiltonian circuit is formed.

**2.3 A Method to Estimate Performance of any Heuristic:** We use any heuristic algorithm to obtain a reasonably good Hamiltonian circuit in a weighted complete graph. We form, *WVAB(G)*, the weighted-vertex-adjacency-bitableau, for the given weighted complete graph. We break the entries in the rows of the right tableau into arrays of columns such that the first array contains labels of vertices containing smallest weight in front of them written in the bracket. The weights of the entries in the successive arrays form the non-decreasing order.

When we get the Hamiltonian tour made up of entries entirely belonging to first array we have assuredly got an exact solution for the TSP. In other words, if the subgraph formed by the edges in the first array is Hamiltonian then every Hamiltonian circuit in it is a shortest Hamiltonian circuit. It is also clear that the value (sum of weights of edges) of any shortest Hamiltonian circuit should be at least equal to the sum of weights of the entries (written in the brackets in front of them) obtained by taking one entry from each row of the first array.

We now systematically describe the procedure of estimating the possible distance of the exact solution and the solution obtained by any heuristic in the following steps:

**Algorithm 2.3.1:**

(1)Construct $WVAB(G)$, which is $VAB(G)$ with weights of the corresponding edges written in the brackets in front of the numbers in the right tableau.

(2)Sort the numbers in the rows of the right tableau and arrange them in non decreasing order of their weights written in the brackets and thus obtain sorted $WVAB(G)$, say $SWVAB(G)$.

(3)Carry out partitioning of $SWVAB(G)$ into arrays of the entries in the right tableau such that the first array contains the entries with smallest weight in all the rows, the second column contains next smallest entries in all the rows, and so on and thus construct table of sorted arrays, say $SWA(G)$

(4)Using any of the heuristics obtain a Hamiltonian circuit (which will contain as many as possible entries belonging to arrays with smaller array numbers, depending upon the performance guarantee of the used heuristic.

(5) Suppose the Hamiltonian circuit formed contains entries from first, second, **...**, $p$-th rows belonging respectively to $i_1$-th , $i_2$-th, **...**, $i_p$-th array and let the difference in weights in the entries in the $i_1$-th , $i_2$-th, **...**, $i_p$-th arrays and first array in the respective first, second, **...**, $p$-th rows be $w_{i_1}, w_{i_2}, \cdots, w_{i_p}$ respectively then the circuit obtained could differ from the shortest circuit by amount **at most** equal to $\sum_{j=1}^{p} w_{i_j}$ .

□

**Example 2.3.1:** Consider the following $WVAB(K_6)$ for the weighted complete graph on six points:

$$WVAB(K_6) = \left[\begin{array}{c|ccccc} 1 & 2(2) & 3(3) & 4(4) & 5(1) & 6(1) \\ 2 & 1(2) & 3(1) & 4(3) & 5(2) & 6(3) \\ 3 & 1(3) & 2(1) & 4(4) & 5(3) & 6(4) \\ 4 & 1(4) & 2(3) & 3(4) & 5(4) & 6(3) \\ 5 & 1(1) & 2(2) & 3(3) & 4(4) & 6(2) \\ 6 & 1(1) & 2(3) & 3(4) & 4(3) & 5(2) \end{array}\right]$$

$$SWVAB(K_6)=\left[\begin{array}{c|ccccc} 1 & 5(1) & 6(1) & 2(2) & 3(3) & 4(4) \\ 2 & 3(1) & 1(2) & 5(2) & 4(3) & 6(3) \\ 3 & 2(1) & 1(3) & 5(3) & 4(4) & 6(4) \\ 4 & 2(3) & 6(3) & 1(4) & 3(4) & 5(4) \\ 5 & 1(1) & 2(2) & 6(2) & 3(3) & 4(4) \\ 6 & 1(1) & 5(2) & 2(3) & 4(3) & 3(4) \end{array}\right]$$

$$SWA(K_6)=\left[\begin{array}{c|cc|cc|cc|c} 1 & 5(1) & 6(1) & 2(2) & & 3(3) & & 4(4) \\ 2 & 3(1) & & 1(2) & 5(2) & 4(3) & 6(3) & \\ 3 & 2(1) & & 1(3) & 5(3) & 4(4) & 6(4) & \\ 4 & 2(3) & 6(3) & 1(4) & 3(4)5(4) & & & \\ 5 & 1(1) & & 2(2) & 6(2) & 3(3) & & 4(4) \\ 6 & 1(1) & & 5(2) & & 2(3) & 4(3) & 3(4) \end{array}\right]$$

Now the sum of weights of entries in the first array, obtained by taking weight of (any) one entry from each row is equal to 8 units, so, the shortest Hamiltonian circuit (when it can be formed using entries in the first array) will have weight at least equal to 8 units. It is easy to see that we can't form a Hamiltonian circuit using entries only from first array. Now, using some approximation algorithm suppose we obtain the circuit as follows:

$1 \rightarrow 6 \rightarrow 5 \rightarrow 2 \rightarrow 3 \rightarrow 4 \rightarrow 1$. Then, this Hamiltonian circuit **could be** away from the exact solution at most by $0+1+1+0+3+1=6$ units.

**2.4 An Exact Algorithm for TSP using *f* –Cutset Matrix:** We construct a lattice, lattice-cutset-graph, and proceed to obtain shortest-extended-induced-cutset-tree, by any efficient shortest tree finding algorithm. We proceed in the following steps:

**Algorithm 2.4.1:**

(1) Construct *f*-cutset matrix,

$$C_f = [C_c : I_{p-1}]_{q\times(p-1)}$$

with respect to some spanning tree for the given weighted complete graph.

(2) Construct lattice, of size $q \times (p-1)$, made up of dots.

(3) Construct lattice-cutset-graph using the given *f*-cutset matrix.

(4) Assign weight equal to the weight of the edge in the originally given weighted complete graph represented by a column in $C_c$ to first (or some) vertical edge belonging to lattice-cutset-graph in that column and assign weight “0” to all other vertical edges in that column (in order to avoid the multiple counting of that weight of that edge represented by the vertex pair written in the top row in the original graph in the desired shortest Hamiltonian circuit) and do this procedure for all columns.

(5) Assign weight equal to the weight corresponding to a branch to all edges reaching a vertex in the last $(p-1)$ columns (branch part) and do the same for all branches.

(7) Assign weight equal to “zero” to all other horizontal edges connecting the points belonging to $C_c$ part when there exists at least one vertical edge emerging from its both ends.

(8) Assign weight to a horizontal edge connecting the points belonging to $C_c$ part equal to the weight of the edge of the chord represented by the column in which the pendant point belongs and has no vertical edge emerging from it.

(9) Choose columns containing smallest nonzero weight for the vertical edge (note that when we select a vertical edge we have by this act actually chosen a vertical path joining all the “1s” in succession on that column). Then at each stage of selection select the next column containing smallest nonzero edge among the remaining columns such that the choice leads to formation of **shortest** induced-cutset-tree.

(10) Find the **shortest** extended-induced-cutset-tree in this graph, by extending the induced-cutset-tree by horizontally going edges to the vertex in the last $(p-1)$ columns in those rows containing odd number of vertices in the first $(q-p+1)$ columns which will lead to formation of a desired Hamiltonian circuit.

**Remark 2.4.1:** Thus, the algorithm essentially consists of finding **shortest-induced-cutset-tree** in the lattice-cutset-graph by any efficient algorithm for finding shortest tree, **similar to** the one due to Kruskal [7] or Prim [8]. Actually this tree is not spanning i.e. does not contain every vertex (vertices in the lattice-cutset-graph are actually representing edges of the original graph present in different cutsets and the number of vertices associated with the same edge is equal to number of *f*-cutsets in which the edge is present) of lattice-cutset-graph. Also, first we have to obtain shortest-induced-cutset-tree, and it is important to note that its determination is **not degree constrained**, and later the extended- shortest-induced-cutset-tree to achieve evenness (the degree constraint comes here but the procedure does not become difficult due to this constraint, here one only need to extend the shortest-induced-cutset-tree by appending the branches wherever required to maintain the required evenness of selected vertices on a row) by adding of selected vertices in each row of lattice graph.

**Remark 2.4.2:** Note that assigning weights to the edges of extended-induced-cutset-tree is set up in a way described in the steps of the algorithm essentially to include the weights of actual chords and branches in the original graph forming the desired Hamiltonian circuit (again, in the original graph) only once.

**Example 2.4.1:** Consider weighted complete graph of **example 2.2.1**, and take any tree, say $T$ , as follows:

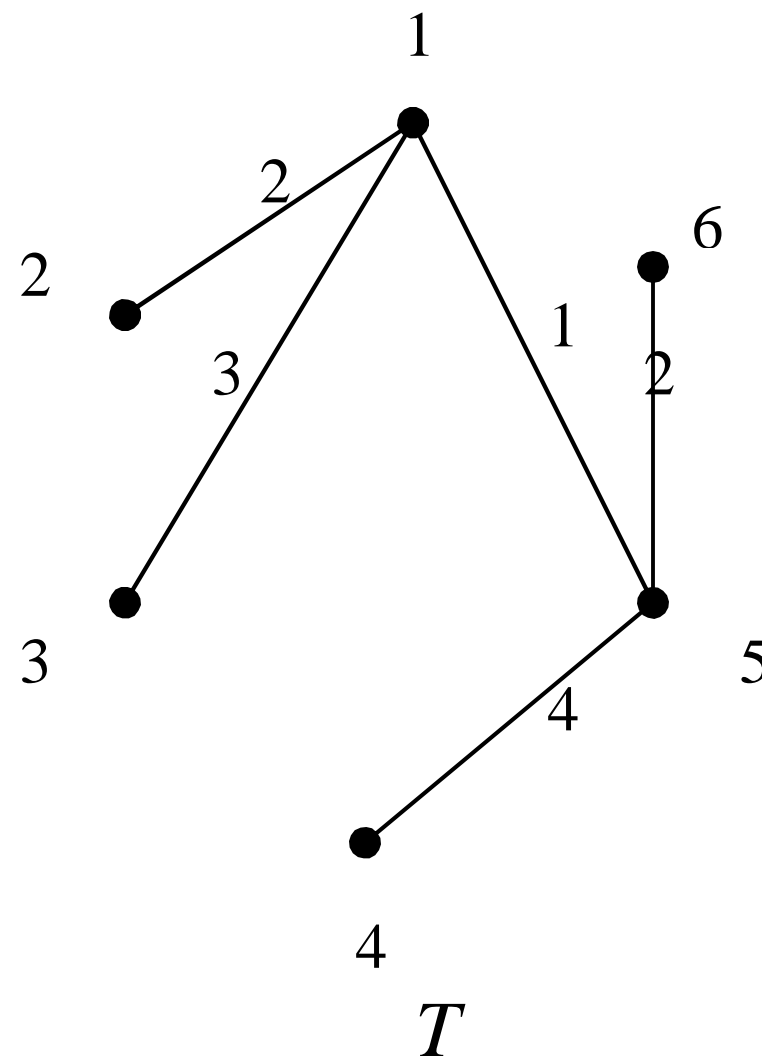

The corresponding *f*-cutset matrix for this tree is given below:

$C_f$ =

| | (1,4) | (1,6) | (2,3) | (2,4) | (2,5) | (2,6) | (3,4) | (3,5) | (3,6) | (4,6) | (1,2) | (1,3) | (1,5) | (4,5) | (5,6) |
|---|---|---|---|---|---|---|---|---|---|---|---|---|---|---|---|
| $C_1$ | 0 | 0 | 1(1) | 1(3) | 1(2) | 1(3) | 0 | 0 | 0 | 0 | 1(2) | 0 | 0 | 0 | 0 |
| $C_2$ | 0 | 0 | 1(1) | 0 | 0 | 0 | 1(4) | 1(3) | 1(4) | 0 | 0 | 1(3) | 0 | 0 | 0 |
| $C_3$ | 1(4) | 1(1) | 0 | 1(3) | 1(2) | 1(3) | 1(4) | 1(3) | 1(4) | 0 | 0 | 0 | 1(1) | 0 | 0 |
| $C_4$ | 1(4) | 0 | 0 | 1(3) | 0 | 0 | 1(4) | 0 | 0 | 1(3) | 0 | 0 | 0 | 1(4) | 0 |
| $C_5$ | 0 | 1(1) | 0 | 0 | 0 | 1(3) | 0 | 0 | 1(4) | 1(3) | 0 | 0 | 0 | 0 | 1(2) |

The first row in the above matrix contains pair of vertices noting edges while the first column records the labels of fundamental cutsets.
As per the above algorithm 2.4.1 columns corresponding to edges {(2,3), (3,4), (4,6)} together form the desired shortest induced-cutset-tree and when appended appropriately by columns corresponding to edges {1,2) (1,5), (1,6)} we form shortest extended-induced-cutset-tree, as shown below:

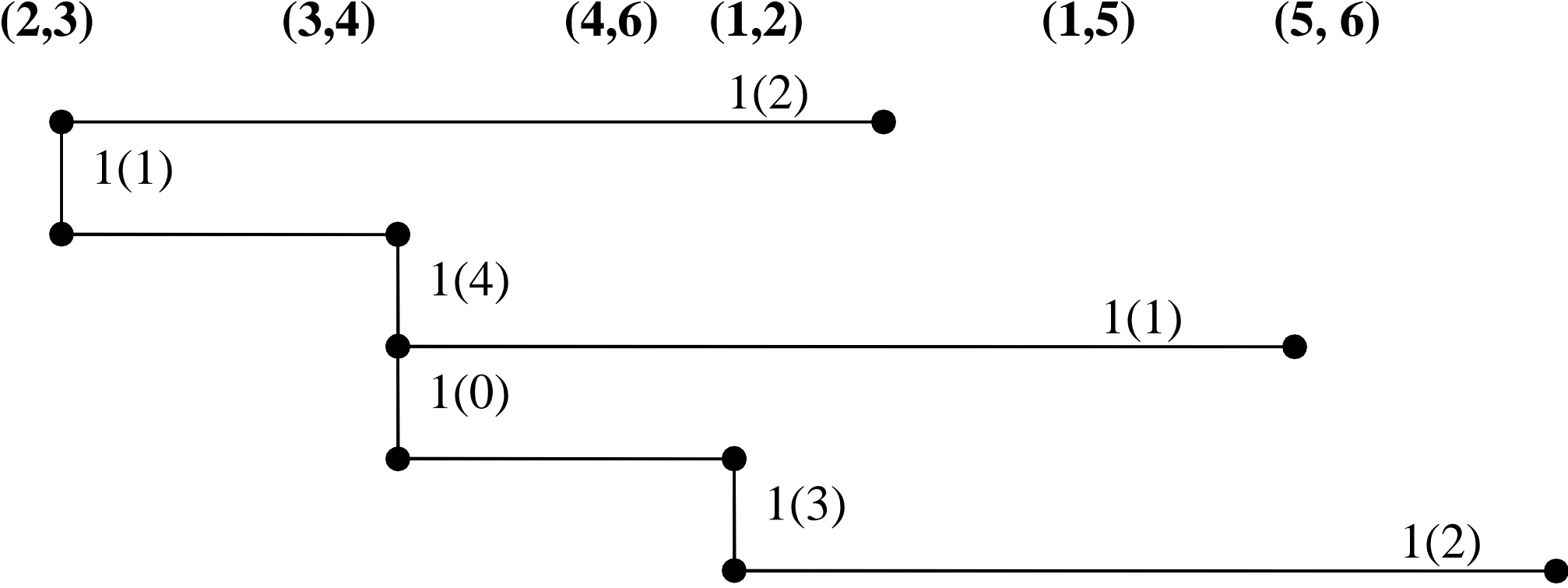

Thus, the Hamiltonian circuit contains edges:
{(2, 3), (3, 4), (4, 6), (1, 2), (1, 5), (5, 6)}
and the actual Hamiltonian circuit obtained by this algorithm is:

2 1 4 3 2 1

$1 \rightarrow 2 \rightarrow 3 \rightarrow 4 \rightarrow 6 \rightarrow 5 \rightarrow 1$, giving total weight "13 Units". This is actually the shortest Hamiltonian circuit in the given complete graph!

**Remark 2.4.2:** From the above algorithm 2.4.1, the problem of finding shortest weight **Hamiltonian circuit** thus reduces to finding **certain shortest tree** (using appropriate columns in the *f*-cutset matrix) in the newly constructed associated **lattice graph** as per this algorithm!

**2.5 A New Polynomial Algorithm for finding Shortest Hamiltonian Path and circuit in a Weighted Complete Graph:**

In this section we begin with proposing two approximation algorithms for shortest Hamiltonian graphs which essentially consists of applying certain chosen permutations (transpositions or product of transpositions) on the adjacency matrix of given weighted complete graph causing reshuffling of the labels of its vertices. We change the labels of vertices through proper choice of permutations in such a way that in this relabeled graph the Hamiltonian path 1→2→3→….k→(k+1)→…→*p* becomes approximation to shortest path in the given weighted complete graph under consideration. We then define so called ordered weighted adjacency list for given weighted complete graph and proceed to the main result of the paper, namely, the exact algorithm based on utilization of ordered weighted adjacency list and the simple properties that any path or circuit must satisfy. This algorithm performs checking of sub-lists, containing (n-1) entries (edge pairs) for paths and n entries (edge pairs) for circuits, chosen from ordered adjacency list in a well defined sequence to determine exactly the shortest Hamiltonian path and shortest Hamiltonian circuit. The procedure has intrinsic advantage of landing on the desired solution in quickest possible time and even in worst case in polynomial time.

Let *G* be a weighted complete graph with the vertex set $V(G)$ and edge set $E(G)$ respectively:

$$V(G) = \{ v_1, v_2, \cdots, v_p \} \text{ and}$$
$$E(G) = \{ e_1, e_2, \cdots, e_q \}$$

Let $A_G = [w_{ij}]_{p\times p}$ denotes the weighted adjacency matrix of *G*.

**Note:** Applying transposition (m, n) on $A_G$ is essentially equivalent to interchanging rows as well as columns, m and n. That is replace m-th row in $A_G$ by n-th row and vice versa and then in thus transformed matrix replace m-th column by n-th column and vice versa (order of these operations, i.e. whether you interchange rows first and then

interchange columns or you interchange columns first and then interchange rows, is immaterial as it produce same end result). Note that this transformation essentially produces a new weighted adjacency matrix that will result due to interchanging labels of vertices $v_m, v_n$ in the original weighted complete graph. We now discuss following algorithm which essentially is an approximation algorithm that produce the result comparable to one obtains from known approximation algorithms.

**Algorithm 2.5.1 (An Approximation Algorithm):**

**(1)** If entry at position (1, 2) in the matrix, i.e. weight $w_{12}$ is already smallest in the first row then proceed to step 2. Else, among the weights $w_{1j}$, $j = 2,3,..., p$ , find minimum weight, say $w_{1j_1}$ . Apply transposition $(2, j_1)$ on $A_G$ , producing new weighted adjacency matrix, say $A_{G_1}$ .

**(2)** If entry at position (2, 3) in the matrix, i.e. weight $w_{23}$ is already smallest in the second row then proceed to step 3. Else, among the weights $w_{2j}$, $j = 3,4,..., p$ , find minimum weight, say $w_{2j_2}$ . Now apply transposition $(3, j_2)$ on $A_{G_1}$ , producing new weighted adjacency matrix, say $A_{G_2}$ .

**(3)** If entry at position (3, 4) in the matrix, i.e. weight $w_{34}$ is already smallest in the third row then proceed to step 4. Else, among the weights $w_{3j}$, $j = 4,5,..., p$ , find minimum weight, say $w_{3j_3}$ . Now apply transposition $(4, j_3)$ on $A_{G_2}$ , producing new weighted adjacency matrix, say $A_{G_3}$ .

**(4)** Continue this procedure applying appropriate transpositions till we finally reach (*p*-2)-th row and among the weights $w_{(p-2)j}$, $j = (p-1), p$ , find minimum weight, say $w_{(p-2)j_{(p-2)}}$ . Now apply transposition $((p-1), j_{(p-2)})$ on $A_{G_{(p-3)}}$ , producing new weighted adjacency matrix, say $A_{G_{(p-2)}}$ .

**(5)** Find the sum of weights of edges in the Hamiltonian path

$$1 \rightarrow 2 \rightarrow 3 \rightarrow \cdots \rightarrow j \rightarrow (j+1) \rightarrow \cdots \rightarrow (p-1) \rightarrow p$$

□

**Remark:** After carrying out "algorithm 1" on given weighted complete graph the Hamiltonian path

$$1 \to 2 \to 3 \to \cdots \to j \to (j+1) \to \cdots \to (p-1) \to p$$

produces a good **approximation** for shortest Hamiltonian path in the given (and conveniently relabeled due to applied transpositions) weighted complete graph. This Hamiltonian path thus obtained will not necessarily be a shortest one.

What we have necessarily achieved is as follows: By application of permutation (transposition) we bring smallest weight entry in the first row at position (1, 2) in the weighted adjacency matrix. This is achieved by transposition of type $(2, j_1)$, where $j_1 > 2$. The algorithm then applies transposition which brings smallest weight entry in the second row at position (2, 3), in the transformed weighted adjacency matrix that results after applying transposition mentioned above. This is achieved by transposition of type $(3, j_2)$, where $j_2 > 3$. Note that because of its special form this second transposition doesn't affect the smallest entry achieved at position (1, 2) while bringing smallest entry (weight) in the second row at position (2, 3) by this second transposition! This story continues, i.e. the later applied transpositions doesn't affect the results of earlier transpositions because of the special choice of the transpositions and at end achieves smallest possible weights in the rows at positions on the **diagonal neighboring the principle diagonal**, i.e. at positions (1, 2), (2, 3), …., (*p*-1, *p*), of the evolved weighted adjacency matrix, evolved through the successive transpositions of specially chosen type. Thus, we have achieved the neighboring diagonal which represents the weights on the Hamiltonian path

$$1 \to 2 \to 3 \to \cdots \to j \to (j+1) \to \cdots \to (p-1) \to p$$

to contain smallest entries from the rows of initially given weighted adjacency matrix.

**Question 2.5.1:** When the Hamiltonian path

$$1 \to 2 \to 3 \to \cdots \to j \to (j+1) \to \cdots \to (p-1) \to p$$

thus produced by this algorithm will be the desired shortest Hamiltonian path?
**Answer:** The Hamiltonian path

$$1 \to 2 \to 3 \to \cdots \to j \to (j+1) \to \cdots \to (p-1) \to p$$

will be shortest if and only if we will (somehow) manage the maximization of sum of weights of entries in the triangular shaped submatrix of the transformed adjacency matrix, i.e. when the following sum

$$\sum_{i=1}^{p-2} \sum_{j=(i+2)}^{p} w_{ij}$$

gets **maximized.**

Now, the next question that naturally arises is as follows:

**Question 2.5.2:** How to maximize this sum?

We will try to come to its possible answer but before that let us consider following

**Example 2.5.1:** We consider following weighted adjacency matrix representing a weighted complete graph and find the Hamiltonian path in its relabeled copy in the form

$$1 \rightarrow 2 \rightarrow 3 \rightarrow \cdots \rightarrow j \rightarrow (j+1) \rightarrow \cdots \rightarrow (p-1) \rightarrow p$$

by applying "algorithm 1". We will see that this Hamiltonian path is not shortest one. We consider the following weighted adjacency matrix and apply "algorithm 1" to it:

$$\begin{bmatrix} 0 & 1 & 6 & 8 & 4 \\ 1 & 0 & 8 & 5 & 6 \\ 6 & 8 & 0 & 9 & 7 \\ 8 & 5 & 9 & 0 & 8 \\ 4 & 6 & 7 & 8 & 0 \end{bmatrix}$$

(1) Since entry at position (1, 2) is already smallest in the first row we proceed to next step.

(2) Since entry in position (2, 4) = 5 is smallest in second row we apply transposition (3, 4) on the above matrix that results into matrix

$$\begin{bmatrix} 0 & 1 & 8 & 6 & 4 \\ 1 & 0 & 5 & 8 & 6 \\ 8 & 5 & 0 & 9 & 8 \\ 6 & 8 & 9 & 0 & 7 \\ 4 & 6 & 8 & 7 & 0 \end{bmatrix}$$

(3) Since entry in position (3, 5) = 8 is smallest in third row we apply transposition (4, 5) on the above matrix that results into matrix

$$\begin{bmatrix} 0 & 1 & 8 & 4 & 6 \\ 1 & 0 & 5 & 6 & 8 \\ 8 & 5 & 0 & 8 & 9 \\ 4 & 6 & 8 & 0 & 7 \\ 6 & 8 & 9 & 7 & 0 \end{bmatrix}$$

Clearly, in this transformed weighted adjacency matrix the Hamiltonian path

$$1 \to 2 \to 3 \to 4 \to 5$$

has total weight $\sum_{i=1}^{4} w_{i,(i+1)} = 21$

Now, it is easy to check that $\sum_{i=1}^{p-2} \sum_{j=(i+2)}^{p} w_{ij}$ in this case is equal to 41.

This sum is actually **not maximized** as we will see below. We actually need to apply some more permutations on the given weighted adjacency matrix that are suitable to maximize this sum. Now without displaying all necessary suitable permutations we have to carry out for maximization of this sum we give the final result below depicting the transformed form of the **same** weighted adjacency matrix with which we started applying algorithm 1. It is

$$\begin{bmatrix} 0 & 5 & 8 & 8 & 9 \\ 5 & 0 & 1 & 6 & 8 \\ 8 & 1 & 0 & 4 & 6 \\ 8 & 6 & 4 & 0 & 7 \\ 9 & 8 & 6 & 7 & 0 \end{bmatrix}$$

Clearly, in this transformed weighted adjacency matrix the Hamiltonian path

$$1 \to 2 \to 3 \to 4 \to 5$$

has total weight $\sum_{i=1}^{4} w_{i,(i+1)} = 17$

Now, it is easy to check that $\sum_{i=1}^{p-2} \sum_{j=(i+2)}^{p} w_{ij}$ in this case is equal to 45.

It can be checked by **brute force** that is this desired sum has **maximized** now and so the Hamiltonian path that we have thus obtained in the transformed weighted adjacency matrix is actually the **shortest** one!

Thus, the problem of finding shortest Hamiltonian path in the form

$$1 \to 2 \to 3 \to \cdots \to j \to (j+1) \to \cdots \to (p-1) \to p$$

has become a **constrained optimization problem** of following type:

**Problem 2.5.1:** Given weighted adjacency matrix corresponding to a given weighted complete graph then devise permutations which will transform this matrix such that path

$$1 \to 2 \to 3 \to \cdots \to j \to (j+1) \to \cdots \to (p-1) \to p$$

has shortest length. In other words, devise sequence of permutations to be applied on given weighted adjacency matrix to be applied on the given weighted complete graph such that $\sum_{i=1}^{p-2} \sum_{j=(i+2)}^{p} w_{ij}$ gets maximized and thus the transformed matrix represents the same weighted complete graph in disguise.

Now, this sum can be seen as made up of sum of entries in columns $p$, ($p$-1), ($p$-2), …, such that $p$-th column contains entries $w_{1p}$, $w_{2p}$, …. , $w_{(p-2)p}$ , ($p$-1)-th column contains entries $w_{1(p-1)}$, $w_{2(p-1)}$, …. , $w_{(p-3)(p-1)}$, etc.

We now proceed to discuss a possible algorithm to tackle this problem.

**Algorithm 2.5.2 (An Approximation Algorithm):**

**(1)** We begin with maximizing $w_{1p}$ . Pick largest weight edge say $w_{ij}$ in the given weighted adjacency matrix, $A_G$ . Transform this weight to position $w_{1p}$ by applying product of transpositions $(1,i)(j,p)$ on $A_G$ . Thus, we have now maximized $w_{1p}$ .

**(2)** Now, among the edges emerging from vertex with label 1 and *p* in the transformed weighted adjacency matrix, $A_G$ , due to step (1), find some edges $(1,i),(p,j),i \neq j$ such that $(1,i)$ is smallest among the edges emerging from 1 and $(p,j)$ is smallest among the edges emerging from *p*. Apply product of transpositions $(2,i)(p-1,j)$ on new transformed $A_G$ we get after step (1) so that weights $w_{12}$ and $w_{(p-1)p}$ are now smallest possible for the situation.

**(3)** Now, among the edges emerging from vertex with label 2 and (*p*-1) in the transformed weighted adjacency matrix, $A_G$ , due to step (2), find some edges $(2,i),(p-1,j),i \neq j$ such that $(2,i)$ is smallest among the edges emerging from 2 and $(p-1,j)$ is smallest among the edges emerging from (*p*-1). Apply product of transpositions $(3,i)(p-2,j)$ on new transformed $A_G$ we get after step (2) so that weights $w_{23}$ and $w_{(p-2)(p-1)}$ are now smallest possible for the situation.

**(4)** Continue steps like (2) and (3) of applying suitable transpositions till we reach at the state having smallest possible weights for $w_{12}$ , $w_{23}$ , $w_{34}$ , …., $w_{(p-2)(p-1)}$ , $w_{(p-1)p}$ in the given situation and all other edges have larger weights.

**(5)** Now, all vertices have been relabeled and we have assigned smallest possible weights to edges comprising the Hamiltonian path

$$1 \rightarrow 2 \rightarrow 3 \rightarrow \cdots \rightarrow j \rightarrow (j+1) \rightarrow \cdots \rightarrow (p-1) \rightarrow p$$

Now, apply suitable transpositions of type $(i, i+2), i = 1,2,\cdots$ by checking that they cause improvement in the sum $\sum_{i=1}^{p-1} w_{i(i+1)}$.

□

**Example 2.5.2:** We consider following weighted adjacency matrix to apply "algorithm 2".

$$\begin{bmatrix} 0 & 11 & 2 & 5 & 3 \\ 11 & 0 & 1 & 6 & 3 \\ 2 & 1 & 0 & 12 & 4 \\ 5 & 6 & 12 & 0 & 8 \\ 3 & 3 & 4 & 8 & 0 \end{bmatrix}$$

**(1)** Since $w_{34}$ is largest so we apply product of transpositions (1, 3)(4,5) on this matrix leading to

$$\begin{bmatrix} 0 & 1 & 2 & 4 & 12 \\ 1 & 0 & 11 & 3 & 6 \\ 2 & 11 & 0 & 3 & 5 \\ 4 & 3 & 3 & 0 & 8 \\ 12 & 6 & 5 & 8 & 0 \end{bmatrix}$$

**(2)** $w_{12}$ ia already minimum so to minimize $w_{45}$ we apply transposition (3, 4) on above matrix leading to new matrix as follows:

$$\begin{bmatrix} 0 & 1 & 4 & 2 & 12 \\ 1 & 0 & 3 & 11 & 6 \\ 4 & 3 & 0 & 3 & 8 \\ 2 & 11 & 3 & 0 & 5 \\ 12 & 6 & 8 & 5 & 0 \end{bmatrix}$$

**(3)** As per step (5) of the algorithm to achieving further minimization, we apply transposition (1, 3) on above matrix leading to new matrix as follows:

$$\begin{bmatrix} 0 & 3 & 4 & 3 & 8 \\ 3 & 0 & 1 & 11 & 6 \\ 4 & 1 & 0 & 2 & 12 \\ 3 & 11 & 2 & 0 & 5 \\ 8 & 6 & 12 & 5 & 0 \end{bmatrix}$$

Clearly, in this transformed weighted adjacency matrix the Hamiltonian path

$$1 \rightarrow 2 \rightarrow 3 \rightarrow 4 \rightarrow 5$$

has total weight $\sum_{i=1}^{4} w_{i,(i+1)} = 11$

**The Ordered Weighted Adjacency List and Sub-lists:** We now proceed with the discussion of the main results of this paper. We propose a smart algorithm to find shortest Hamiltonian path and shortest Hamiltonian circuit in the given weighted complete graph.

**Definition 2.5.1:** The **weighted adjacency list**, $WAL(G)$, associated with the given weighted complete graph, $G$ , on *p* vertices, is the following bitableau in which the left tableau represents weights of the edges represented by vertex pairs written in the same row in the right tableau.

**Definition 2.5.2:** The weighted adjacency list is called **ordered weighted adjacency list** and denoted as $OWAL(G)$ when the rows of weighted adjacency list are so permuted that the weights in the left tableau get ordered, i.e. these weights form a nondecreasing sequence in the downward direction. In other words, the weighted adjacency list becomes ordered weighted adjacency list when the left tableau becomes a nondecreasing column of entries representing weights of the edges written in the corresponding row in the right tableau. The $WAL(G)$ is following bitableau:

$$WAL(G) = \begin{bmatrix} w_{i_1 j_1} & (i_1, j_1) \\ w_{i_2 j_2} & (i_2, j_2) \\ \bullet & \bullet \\ \bullet & \bullet \\ w_{i_l j_l} & (i_l, j_l) \\ \bullet & \bullet \\ w_{i_m j_m} & (i_m, j_m) \end{bmatrix}$$

where, $m = \frac{p(p-1)}{2}$.

The $OWAL(G)$ is the following bitableau:

$$OWAL(G) = \begin{bmatrix} w_{i_1 j_1} & (i_1, j_1) \\ w_{i_2 j_2} & (i_2, j_2) \\ \bullet & \bullet \\ \bullet & \bullet \\ w_{i_l j_l} & (i_l, j_l) \\ \bullet & \bullet \\ w_{i_m j_m} & (i_m, j_m) \end{bmatrix}$$

where $m = \frac{p(p-1)}{2}$ and in addition, $w_{i_1 j_1} \leq w_{i_2 j_2} \leq \cdots \leq w_{i_m j_m}$

**Definition 2.5.3:** The **weighted adjacency sub-list**, $SubWAL(G)$, is essentially any sub-bitableau made by picking some portion of the $WAL(G)$, i.e. made by picking any rows among the rows in $WAL(G)$.

**Definition 2.5.4:** The **ordered weighted adjacency sub-list**, $SubOWAL(G)$, is essentially any sub-bitableau made by picking some portion of the $OWAL(G)$ and keeping them in the same nondecreasing order, i.e. made by picking any rows among the rows in $OWAL(G)$ and keeping them in the same ordered form.

**Definition 2.5.5:** The **size** of weighted adjacency sub-list, or ordered weighted adjacency sub-list, is the cardinality of entries in the sub-list, i.e. number of rows in the sub-bitableau representing that sub-list.

**Definition 2.5.6:** The **weight** of the weighted adjacency sub-list, or ordered weighted adjacency sub-list, is the sum of weights in the left sub-tableau representing that sub-list.

It is **easy to check** that

(A) A set of ($p$-1) vertex pairs in the right tableau of $OWAL(G)$ together represents a Hamiltonian path if (i) these pairs together contain all the vertices, (ii) the degrees of all but two vertices is two, (iii) the degree of the left out two vertices is one, and (iv) these vertex pairs together form a connected graph.

(B) A set of $p$ vertex pairs in the right tableau of $OWAL(G)$ together represents a Hamiltonian circuit if (i) these pairs together contain all the vertices, (ii) the degrees of all vertices is two, and (iii) these vertex pairs together form a connected graph.

These two easy checks will form an **important backbone** of our exact algorithm for finding shortest Hamiltonian path or shortest Hamiltonian circuit in the given weighted complete graph.

We now proceed with

**Algorithm 2.5.3 Shortest Hamiltonian Path (Exact):**

**(1)** Form ordered weighted adjacency list, $OWAL(G)$, corresponding to given weighted complete graph.

**(2)** Form all possible ordered weighted adjacency sub-lists, $SubOWAL(G)$, each of size ($p$-1).

**(3)** Arrange these sub-lists in lexicographic order in accordance with their respective weights.

**(4)** Use easy check (A) mentioned above in succession (starting with smallest weight sub-list) on each sub-list and stop at the first success.

**(5)** Record the Hamiltonian path thus obtained and its weight. This will be a desired **shortest Hamiltonian path!**

□

**Algorithm 2.5.4 Shortest Hamiltonian Circuit (Exact):**

**(1)** Form ordered weighted adjacency list, $OWAL(G)$, corresponding to given weighted complete graph.

**(2)** Form all possible ordered weighted adjacency sub-lists, $SubOWAL(G)$, each of size $p$.

**(3)** Arrange these sub-lists in lexicographic order in accordance with their respective weights.

**(4)** Use easy check (B) mentioned above in succession (starting with smallest weight sub-list) on each sub-list and stop at the first success.

**(5)** Record the Hamiltonian circuit thus obtained and its weight. This will be a desired **shortest Hamiltonian circuit!**

□

**Example 3:** Consider following weighted adjacency matrix:

$$\begin{bmatrix} 0 & 1 & 6 & 8 & 4 \\ 1 & 0 & 8 & 5 & 6 \\ 6 & 8 & 0 & 9 & 7 \\ 8 & 5 & 9 & 0 & 8 \\ 4 & 6 & 7 & 8 & 0 \end{bmatrix}$$

The ordered weighted adjacency list for this is as below, conveniently in the form of a 2-columned table:

| | |
|---|---|
| 1 | (1,2) |
| 4 | (1,5) |
| 5 | (2,4) |
| 6 | (1,3) |
| 6 | (2,5) |
| 7 | (3,5) |
| 8 | (2,3) |
| 8 | (1,4) |
| 8 | (4,5) |
| 9 | (3,4) |

It is easy to check that pairs {(1,2), (1,5), (2,4), (3,5)} together form the desired shortest Hamiltonian path as per the "algorithm 3" and pairs {(1,2), (1,5), (2,4), (3,5), (3,4)} together form the desired shortest Hamiltonian circuit as per the "algorithm 4".

**2.6 A New TSP Algorithm for Graphs satisfying Triangle Inequality:** For complete graphs with vertices as points representing cities on a realistic planar map and where the weights on the edges are real distances we propose the following criterion for a Hamiltonian tour to be optimal. This criterion is motivated from appropriate extension of the idea of geodesic, i.e. the curve of shortest length joining two points on a surface that we determine for the surface under consideration using calculus of variation, e.g. by solving Euler-Lagrange equation one finds that straight line is geodesic for plane and for any two points on plane, the length of the straight line segment joining these two points is thus the shortest distance.

Now, suppose we are given three non-collinear points on the plane say A, B, C. What is the smooth curve of shortest length (geodesic) passing through all these points? Let us join points A, B and then B, C by a straight line segments (and make smooth their join at B) and obtain smooth curve A-B-C. Similarly, join points A, C and then C, B by a straight line segments (and make smooth their join at C) and obtain smooth curve A-C-B. Now which one among these curves A-B-C and A-C-B is shortest? Move the tangent vector along both the curves and note down the magnitude of rotation (at B and C respectively) for both these curves. **The curve for which the magnitude of the angle of rotation of tangent vector (at B or C) will be least will be the shortest one!**

Given $p$ distinct points on a plane, $\{ v_1, v_2, \cdots, v_p \}$ and suppose we construct all possible Hamiltonian paths between $v_1$ to $v_p$ (passing through all other points only once). Which one among these is shortest? **The one for which the sum of the magnitudes of the angles of rotation of tangent vector (at each intermediate node) will be least will be the shortest Hamiltonian path!** Also, suppose we construct all possible Hamiltonian tours (passing through all points only once). Which one among these is shortest? **The one for which the sum of the magnitudes of the angles of rotation of tangent vector (at each intermediate node) will be least will be the shortest Hamiltonian tour!**

**Algorithm 2.6.1 (An Economical Hamiltonian Tour for Real TS passing through Real Cities):**

Let $\{ r_1, r_2, \cdots, r_i, \cdots, r_p \}$ be the position vectors of points representing cities on plane (with respect to some origin).

(1) Find the following vector (center of mass), $R = \left( \frac{1}{p} \right) \sum_{i=1}^{p} r_i$

(2) Find $r = \max | R - r_i |$

(3) By taking the point representing vector $R$ (center of mass) as center draw a circle with radius equal to $r$ (clearly, all points represented by position vectors $\{ r_1, r_2, \cdots, r_i, \cdots, r_p \}$ will lie inside this circle).

(4) Draw radii to this circle passing through each point representing a position vector $\{ r_1, r_2, \cdots, r_i, \cdots, r_p \}$.

(5) Name the end points of these radii on the circle by $r_i^{/}$ if the radius is passing through $r_i$, for all $i = 1, 2, \cdots, p$.

(6) Form the sequence of points by starting with some end point of some radius and moving to the next end point of radius in succession anticlockwise (or, clockwise) till you reach to the starting point. Let it be $r_{i_1}^{/} \to r_{i_2}^{/} \to \cdots \to r_{i_1}^{/}$

(7) Join (by vectors) and construct Hamiltonian circuit by replacing $r_{i_j}^{/}$ by $r_{i_j}$ in the above sequence to form $r_{i_1} \to r_{i_2} \to \cdots \to r_{i_1}$.

**Remark 2.6.1:** The Hamiltonian circuit thus obtained will have very much close to the minimum the sum of the magnitudes of the angles of rotations of the tangent vector going around this Hamiltonian circuit.

3. **NP-Complete Problems and Unstructured Search:** It is clear to see that we will be able to solve NP-Complete problems if we can devise some novel algorithm for unstructured search which will have superiority of exponential order over the existing algorithms for unstructured search. For search data of size *N*, Grover's [9] quantum search algorithm $\sim O(\sqrt{N})$ certainly speeds up the well known classical algorithm $\sim O(N)$, but only polynomially.

Grover's quantum search algorithm can be adapted quite readily to solve NP-Complete problems albeit again in exponential time, but (only) with reduced exponent compared to what is known classically [10].In this section we will propose a novel classical algorithm which will show improvement of exponential order over existing algorithms. With this algorithm we will certainly be in the position to say that we indeed are in possession of solution to NP-complete problems! The unstructured search problem asks for search of some predefined number, called target, from given unstructured list of numbers. In this paper we propose a novel classical algorithm with complexity $\sim O(LogN)$ for searching the target from unstructured list of numbers. We thus propose a new algorithm, which achieves improvement of exponential order over existing algorithms. Suppose *N* is the largest number in the list then we consider *N* dimensional vector space with Euclidean basis. With each of the numbers in the given unstructured list we associate the unique basis vector among the vectors that form together the Euclidean basis. For example suppose *j* is a number in the list then we associate with this number *j* the unique basis vector in the above mentioned *N*-dimensional vector space, $| j >$ = transpose(0, 0, 0, … , 0, 0, 1, 0, 0, … , 0, 0, 0), where there is entry 1 only at *j*-th place and every where else there is entry 0. We then divide the given list of numbers in two roughly equal parts (i.e. we divide the given bag, B say, containing scrambled numbers in two roughly equal parts and put them in two separate bags, bag B1 and bag B2. We represent the list of numbers in bag B1, bag B2 in the form of equally weighted superposition of basis vectors associated with the numbers contained in these bags, namely, we represent list in Bag B1 (Bag B2) as a single state formed by equally weighted superposition of orthonormal basis states forming Euclidean basis corresponding to numbers in the bag B1 (bag B2), namely, $|\psi_1>$ ($|\psi_2>$). Let *t* be the target number. It will be represented as basis state $|t>$ called target state. We then find the value of scalar product of target state $|t>$ with $|\psi_1>$ (or $|\psi_2>$). It will revel us whether *t* belongs to bag B1 (or bag B2) which essentially enables us to carry out the binary search and to achieve above mentioned $\sim O(LogN)$ complexity!

**3.1 Preamble:** The unstructured search problem requires us to find a particular target item amongst a set of *N* candidates. We can label these *N* candidates by indices $x$ in the range $1 \le x \le N$, and we are supposed to find the index of the sought after target item, $x = t$, say. Now, suppose these $N$ numbers, as tags associated with items as identifiers for these items, be mixed randomly among each other in a bag. Your task is to pick out the target, $x = t$, in fewest possible trials from this bag containing this randomly done mixture of numbers. It is this random mixing of numbers which makes this list of numbers unstructured. This is essentially the well-known so called problem of unstructured search. Also, it is well known that the existing classical algorithm for the solution of this problem has complexity $\sim O(N)$ and Grover's quantum algorithm developed by Lov Grover [9] has complexity $\sim O(\sqrt{N})$.

We now propose a novel classical algorithm which achieves exponential speedup over existing algorithms. The so called binary search algorithm can determine the number which is declared as target from the given list of $N$ numbers with complexity of the order of $\sim O(LogN)$ but this algorithm works only when given list of numbers is sorted (ordered).

For our new algorithm we begin with associating a unique unit vector from the standard Euclidean basis with each number in the given scrambled bag of numbers, thus, we consider each number in the bag as a unique basis vector in the $N$-dimensional Euclidean space when the largest number in the given scrambled list of numbers in the bag is $N$. Thus, each number $j$ in the given list is represented as state $| j >$, a column vector, where,

$$| j >= \begin{pmatrix} 0 \\ 0 \\ 0 \\ \vdots \\ 0 \\ 1 \\ 0 \\ \vdots \\ 0 \end{pmatrix}$$

It is important to note that all the scalar components of this vector are 0 except scalar component in the place $j$ from top which is equal to 1. Further, we associate equally weighted superposition of basis vectors with given list of numbers. Thus, we represent entire list as a single superposed state. Let $i_1, i_2, i_3, \cdots, i_R, \cdots, i_N$ be the numbers in the bag and there is no order relation among these $N$ numbers, then we represent this entire list as a single state, $| \psi >$, obtained as superposition state:

$$| \psi >=| i_1 > + | i_2 > + | i_3 > + \cdots + | i_R > + \cdots + | i_N >$$

We further consider projection of state $| \psi >$ on the target state, $| t >$. In other words, we consider scalar product $< t | \psi >$. This scalar product will obviously satisfy the following conditions:

$< t \mid \psi > \; = \; 0$, if state $\mid t >$ doesn't belong to basis states forming state $\mid \psi >$

and

$< t \mid \psi > \; = \; 1$, if state $\mid t >$ belongs to basis states forming state $\mid \psi >$

This simple fact will enable us to carry out the desired binary search on the scrambled list of numbers. We will see how in the algorithm given below.

## 3.2 Algorithm:

1. Let B be the bag containing scrambled list of numbers $\{ i_1, i_2, i_3, \cdots, i_R, \cdots, i_N \}$ and let *N* be the largest among the numbers in this given bag B, and let *t* be the number which is target to be searched.
2. Consider *N*-dimensional Euclidean vector space, E, with Euclidean basis of unit vectors. Thus, a basis state $\mid j > \; = \;$ transpose (0, 0, 0, …, 0, 1, 0, …, 0), where all components are 0 except *j*-th component which is 1. Note that with each number *j* in the bag we will be associating state $\mid j >$ in the Euclidean basis.
3. Divide list in bag B, in any arbitrary way, into two sub-lists of roughly equal sizes $N_1$ and $N_2$, and put them in bags B1 and B2. Let now bag B1 contain scrambled numbers $\{ i_{11}, i_{12}, \cdots, i_{1p} \}$ and bag B2 contains scrambled numbers $\{ i_{21}, i_{22}, \cdots, i_{2q} \}$, where, $p \approx q$ from the original scrambled list.
4. Form any one of the following superposed states:

$$\mid \psi_1 >= \mid i_{11} > + \mid i_{12} > + \cdots + \mid i_{1p} >,$$

$$\text{(or, } \mid \psi_2 >= \mid i_{21} > + \mid i_{22} > + \cdots + \mid i_{2q} > \text{)}$$

   and consider scalar product $< t \mid \psi_1 >$ (or $< t \mid \psi_2 >$)
5. Suppose we have formed state $\mid \psi_1 >$ and if $< t \mid \psi_1 > \; = \; 1$ then state $\mid t >$ is member of the superposition of states that forms state $\mid \psi_1 >$. Else, if $< t \mid \psi_1 > \; =$ 0 then state $\mid t >$ is member of the superposition of states that forms state $\mid \psi_2 >$.
6. When $< t \mid \psi_1 > \; = \; 1$, then set B1 = B and go to step 3. When $< t \mid \psi_1 > \; = \; 0$, then set B2 = B and go to step 3.
7. Continue till We (obviously) reach the desired target state and thus the number which is target.

3.3 **Example:** Suppose we are given {2, 11, 7, 5, 3, 6, 9, 4} as scrambled list of numbers in Bag B and suppose number 3 is our target, i.e. we wish to locate and find number 3 in this scrambled list. As per steps of algorithm we divide these numbers, in bag B say, into two bags B1, B2 such that bag B1 contains numbers {2, 11, 7, 5} and bag B2 contains numbers {3, 6, 9, 4}. Since 11 is the largest number in the list of given numbers we consider 11-dimensional Euclidean vector space and form state $|\psi_1>=|2>+|11>+|7>+|5>$, Further, we find scalar product $<3|\psi_1>$, which is equal to 0. So, clearly, target 3 belongs to bag B2. So, we set B2 = B and proceed with step 3, i.e. the division of this newly defined bag B. In other words, we will then have newly defined bags B1 = {3, 6} and B2 = {9, 4}. We form state $|\psi_1>=|3>+|6>$ and again find scalar product $<3|\psi_1>$ which is equal to 1. So, clearly, target 3 belongs to bag B1. So, we set B1 = B and proceed with division of this newly defined bag B. In other words, we will have now newly defined bags B1 = {3} and B2 = {6}. We form state $|\psi_1>=|3>$ and again find scalar product $<3|\psi_1>$ which turns out to be equal to 1. So, we have thus obtained the desired target!

### 3.3.1 Some New Quantum Algorithms for Unstructured search:

We develop three new quantum algorithms for searching the desired target state in the unstructured database of size *N*. The first algorithm requires *Log N* iterative steps. It constructs two quantum bags of equal size in terms of two quantum states, out of which exactly one quantum state will have nonzero overlap with the target state. This determination of overlap is done by taking the inner product, in *Log N* time [11], of the implicitly known target state with any one of these two quantum states. The second algorithm requires just one single step which uses a new suitable operator and the choice of this operator is problem dependent, i.e. it depends upon the number of qubits required to be used to represent an element in the index set. The third algorithm again requires only a single step and this algorithm makes use of a fixed (same) operator. It is known that algorithm for unstructured database search can be easily adaptable for solving *NP-Complete* problems. However, the computational complexity of *NP-Complete* problems after the adaptations of both the classical as well as quantum [9] search algorithms remains of the exponential order as the exponent for quantum [9] algorithm changes only to one-half times the exponent for classical algorithm. But for our quantum algorithms the exponent falls substantially so that our new quantum algorithms for unstructured search are capable if reducing the computational complexity of *NP-Complete* problems to polynomial order! We propose in all three new quantum algorithms for unstructured database search. If *N* is the size of the unstructured database then we show that we can pick out the desired target in just $LogN$ steps by the first algorithm, and in just single step by the other two algorithms! The innovation in the first algorithm consists of dividing the given database into two equal sized databases in terms of two quantum states and By using the idea of taking inner product in $LogN$ time [11] of any one these quantum states representing the quantum bags with the target state which enables one to find out the quantum state to which the target state belongs!

We now proceed to propose our first quantum algorithm for unstructured search. This new quantum algorithm proceeds roughly as follows: It begins with the preparation of the implicitly known desired target state. Starting with a quantum bag that contains target, i.e. starting with a quantum state that contains the target state it then carries out the construction of two suitable initial quantum states using state that contains the target state in the superposition. It then evaluates the inner product of the target state with any one of the two initially constructed quantum states mentioned above. The value of this inner product determines to which quantum state the target state belongs. This quantum state to which target state belong is used further to construct two more new suitable quantum states and the same procedure is repeated iteratively. By iterating these steps for $LogN$ times we will see that with these steps one directly arrives at the desired target state and completes the search. As mentioned above the generate-and-test type classical algorithm or quantum [9] algorithm for unstructured database search though can be easily adapted to solve the *NP-Complete* problems still the computational complexity of these algorithms after the adaptations remains that of the exponential order as the exponent for quantum [9] algorithm changes only to one-half times the exponent for classical algorithm. But for our first quantum algorithm the exponent becomes the polynomial of logarithm of the exponent for the classical algorithm. Therefore, our first quantum algorithm reduces the computational complexity of *NP-Complete* problems to polynomial order! If $N$ is the size of the unstructured database then we show that we can attain the desired target in just $LogN$ steps! To attain the desired target the best known generate-and-test type classical algorithm and quantum [9] algorithm for unstructured search requires roughly $\frac{N}{2}$ steps and $\sqrt{N}$ steps respectively. This implies that only quadratic speedup is achievable by quantum [9] algorithm over classical algorithm. Though such speedup is quite good one still it is not good enough as it doesn't tame the problems with exponential complexity. A formal statement of unstructured search problem is as follows: Consider a search problem that requires to find a particular element of the database. Given a set containing $N$ candidates, and suppose these $N$ candidates are labeled by indices, $x$ in the range $0 \leq x \leq (N-1)$, and that the index of the sought after target item is $x = t$ . Let there be a computational oracle, or "black-box function", $f_t(x)$, that when presented with an index $x$ can pronounce on whether or not it is the index of the target. Specifically, $f_t(x)$ is defined such that $f_t(x) = 1$ if $x = t$ and $f_t(x) = 0$ if $x \neq t$ where 1 stands for YES and 0 stands for NO. The search problem is unstructured because there is no discernible pattern to the values of $f_t(x)$ to provide any guidance in finding $x = t$ . Our job is to find index $x = t$ , using fewest calls to the oracle $f_t(x)$. Oracle is nothing but a factitious mathematical device that allows one to estimate the computational cost of some algorithm measured in the units of the "number of calls to this oracle", required to reach the solution. "Oracle" or "black-box function" or "knowledge holder" are synonyms, and if we consider for example the problem of finding name given telephone number what is the oracle? The 'oracle' in this case is the 'telephone directory' itself. We now express the search problem in quantum mechanical language. A quantum analog of the bag of indices can be regarded as an equally weighted superposition of all the indices in the range $0 \leq x \leq (N-1)$, i.e. the quantum state $|\Psi>=\frac{1}{\sqrt{N}}\sum_{x=0}^{N-1}|x>$. Thus, the bag of all indices can be looked upon as a wave function $|\Psi>$ given above. Let us suppose that $N=2^n$. Therefore, using binary representation for all the indices in the bag we can express the wave function representing bag of indices as

$$|\Psi>=\frac{1}{\sqrt{N}}\sum_{i_1,i_2,\cdots i_n=0}^{1}|i_1 i_2\cdots i_n>$$

To prepare such state is in fact a very easy task. For this one just need to take as a starting state a tensor product of *n* number of zero kets, $|0>=\begin{pmatrix}1\\0\end{pmatrix}$, and then apply Hadamard operator, $H$, on each zero ket, $|0>$, in the tensor product. Thus,

$$|\Psi>=\frac{1}{\sqrt{N}}\sum_{i_1,i_2,\cdots i_n=0}^{1}|i_1 i_2\cdots i_n>=H^{\otimes n}|00\cdots 0>$$

The implicitly known target state $|t>=|t_1 t_2 t_3\cdots t_n>$, where each $t_i\in\{0,1\}$, can be prepared using the oracle, $f_t(x)$, which gives value 1 when $x=t$ and 0 when $x\neq t$, by expressing the target state, $|t>$ corresponding to index $x=t$ by using the relation of the target state, $|t>$, and the oracle function, $f_t(x)$. This relation can be expressed in the following two equivalent forms:

$$|t>=\frac{\sqrt{N}}{2}\left[1-(-1)^{f_t(x)}\right]|\Psi> \qquad \text{(A)}$$

or,

$$|t>=\frac{\sqrt{N}}{2}\left[1-\left[1-2|t><t|\right]\right]|\Psi> \qquad \text{(B)}$$

We now divide the elements in the bag containing indices $x$, such that $0\leq x\leq(N-1)$, into two bags such that the first bag will contain half indices, i.e. all those indices, $x$, such that $0\leq x\leq(\frac{N}{2}-1)$ and the second bag will contain remaining half indices, i.e. all those indices, $x$, such that $\frac{N}{2}\leq x\leq(N-1)$. It is easy to achieve this by constructing these bags in terms of two quantum states, $|\Psi_0>$ and $|\Psi_1>$ as follows, where

$$|\Psi_0>=|0>\otimes H^{\otimes(n-1)}|00\cdots 0>$$

and,

$$|\Psi_1>=|1>\otimes H^{\otimes(n-1)}|00\cdots0>$$

Note that the ket $|00\cdots0>$ in the above expressions for $|\Psi_0>$ and $|\Psi_1>$ is of length $(n-1)$, i.e. a computational basis state in $(2^{(n-1)})$ dimensional Hilbert space, while $|\Psi_0>$ and $|\Psi_1>$ are obviously states in $2^n$ dimensional Hilbert space. Also, $|\Psi_0>$ represents the bag that contains all those indices, $x$, such that $0\leq x\leq(\frac{N}{2}-1)$ and $|\Psi_1>$ represents the bag that contains all those indices, $x$, such that $\frac{N}{2}\leq x\leq(N-1)$, as desired.

The idea behind our new quantum algorithm in simple terms is to divide "the bag which contains the target state" at each iterative step into two separate bags of equal size such that now the target state will belong to some one and only one of these two bags which now has become equal to half of the size of the original bag and then to determine by taking inner product of any one state representing these bags with the target state to which the target state belongs. Thus we manage to reduce the size of the bag that contains the target state in each of the iterations to half of its size at that stage. By proceeding along these lines finally the bag that contains the target state will become of size one, i.e. it will contain only the target state itself. Thus, we first begin with the bag represented by the wave function, $|\Psi>$, this original bag contains all numbers from 0 to *N*-1., i.e. it contains the target state.

We now proceed systematically with our first new quantum algorithm through precise steps as follows:

### 3.4 The First New Quantum Search Algorithm implying *P* = *NP* :

**(i)** Construct quantum state, $|\Psi>$ say, representing the bag of all indices $x$, such that $0\leq x\leq(N-1)$. Namely,

$$|\Psi>=\frac{1}{\sqrt{N}}\sum_{i_1,i_2,\cdots i_n=0}^{1}|i_1i_2\cdots i_n>=H^{\otimes n}|00\cdots0>$$

Let $N=2^n$.Since, $t\in[0,(2^n-1)]$, therefore $|t>\in|\Psi>$, i.e. certainly, $<t|\Psi>\neq0$.

**(ii)** Since $t\in[0,(2^n-1)]$, we divide the indices in this bag into two parts of identical size and put them into two new bags. This is done in equivalent terms as follows. Construct two quantum states $|\Psi_0>$ and $|\Psi_1>$representing these two bags such that $|\Psi_0>$ will represent the bag that contains all those indices, $x$, such that $0 \le x \le (2^{(n-1)} - 1)$ and $|\Psi_1>$ will represent the bag that contains all those indices, $x$, such that $2^{(n-1)} \le x \le (2^n - 1)$. Thus, we have

$$|\Psi_0>=|0>\otimes H^{\otimes(n-1)}|00\cdots0>$$

and,

$$|\Psi_1>=|1>\otimes H^{\otimes(n-1)}|00\cdots0>$$

This further implies that either $|t>\in|\Psi_0>$ or $|t>\in|\Psi_1>$.

**(iii)** Take inner product of the implicitly known target state $|t>$, expressed above in two equivalent forms, (A) or (B), with any one of the two quantum states representing two bags of indices, namely, $|\Psi_0>$ and $|\Psi_1>$ given above.

Case (a): Without loss of generality (WLOG), suppose if $<t|\Psi_0>\neq 0$ then clearly we can infer that $|t>\in|\Psi_0>$, i.e. $t\in[0,(2^{(n-1)}-1)]$, i.e. $t$ belongs to the first bag that contains all those indices, $x$, such that $0 \le x \le (2^{(n-1)} - 1)$.

Case (b): Without loss of generality (WLOG), suppose if $<t|\Psi_0>=0$, i.e. $<t|\Psi_1>\neq 0$, then clearly we can infer that $|t>\in|\Psi_1>$, i.e. $t\in[2^{(n-1)},(2^n-1)]$, i.e. $t$ belongs to the second bag that contains all those indices, $x$, such that $2^{(n-1)} \le x \le (2^n - 1)$.

**(iv)** Case (a): Since $t\in[0,(2^{(n-1)}-1)]$, we divide this bag of indices into two equal parts and put them into two new bags. This is done in equivalent terms as follows. Construct two quantum states $|\Psi_{00}>$ and $|\Psi_{01}>$ representing these two bags such that $|\Psi_{00}>$ will represent the bag that contains all those indices, $x$, such that $0 \le x \le (2^{(n-2)} - 1)$ and $|\Psi_{01}>$ will represent the bag that contains all those indices, $x$, such that $2^{(n-2)} \le x \le (2^{(n-1)} - 1)$. Thus, we have

$$|\Psi_{00}>=|0>|0>\otimes H^{\otimes(n-2)}|00\cdots0>$$

and,

$$|\Psi_{01}>=|0>|1>\otimes H^{\otimes(n-2)}|00\cdots0>$$

This further implies that either $|t>\in|\Psi_{00}>$ or $|t>\in|\Psi_{01}>$.

Case (b): Since $t\in[2^{(n-1)},2^n-1)]$, we divide this bag of indices into two equal parts and put them into two new bags. This is done in equivalent terms as follows. Construct two quantum states $|\Psi_{10}>$ and $|\Psi_{11}>$ representing these two bags such that $|\Psi_{10}>$ will represent the bag that contains all those indices, $x$, such that $2^{(n-1)} \le x \le (2^{(n-1)} + (2^{(n-2)} - 1)$ and $|\Psi_{11}>$ will represent the bag that contains all those indices, $x$, such that $(2^{(n-1)} + 2^{(n-2)}) \leq x \leq (2^n - 1)$. Thus, we have

$$|\Psi_{10}>=|1>|0> \otimes H^{\otimes(n-2)}|00\cdots 0>$$

and,

$$|\Psi_{11}>=|1>|1> \otimes H^{\otimes(n-2)}|00\cdots 0>$$

This further implies that either $|t> \in |\Psi_{10}>$ or $|t> \in |\Psi_{11}>$.

**(v)** As is done in **(iii),** by taking inner product of the target state $|t>$ now with $|\Psi_{00}>$ or $|\Psi_{01}>$ when case (a) is true, or with $|\Psi_{10}>$ or $|\Psi_{11}>$ when case (b) is true we determine to which quantum bag represented by these quantum states the target state is part of, i.e. the target state has a nonzero overlap with. We continue on these lines with dividing, each time the correct quantum bag (the one containing the target state), into two separate new quantum bags till (assuredly) the size of the correct quantum bag (that has nonzero overlap with target state) will reduce to the bag containing just one entry, i.e. the target state itself!!

□

### 3.5 The Second New Quantum Search Algorithm implying *P = NP* :

Let $x$, $0 \leq x \leq (N-1)$, be an element in the unstructured database of size $N$. Let $N = 2^n$ hence $0 \leq x \leq (2^n - 1)$. Our aim in the unstructured database search problem is to locate and pick out the target index, $t \in [0,(2^n - 1)]$. Note that with each index $x$ we can associate a computational basis state, $|x>$ made up of $n$ qubits, i.e. $|x>=|x_1 x_2 \cdots x_n>$, where, $x_i \in \{0,1\}, 1 \leq i \leq n$. So, our aim in the unstructured database search problem is to substantially amplify the amplitude of the target state, $|t>$. Suppose we have a 1-YES quantum oracle defined in terms of operator, $O$, which performs the operation $O|x>=(-1)^{f_t(x)}|x>$, where as mentioned previously $f_t(x)$ is defined such that $f_t(x)=1$ if $x=t$ and $f_t(x)=0$ if $x \neq t$. where $f_t(x)=1$ stands for YES and $f_t(x)=0$ stands for NO. It is clear to check that the operator $O$ is unitary. We can see that the real operator $O$ is an inversion operator which only changes the sign of the target state $|t>$ and keeps all other states $|x>$ unchanged. If we take a wave function, $|\Psi>$ say, made up of some superposition of computational basis states and operate the operator $O$ on it then by its definition it will leave all the computational basis states as they are and will change the sign only that of the computational basis state which is the target state. Now if we will operate $O$ one more time then again it will leave all the computational basis states as they are and will restore the sign of the target state. Thus, $O^2 = OO^+ = O^+O = I$ .We define $M_k = \left[(2^k |\Psi><\Psi| - (2^k - 1)I)\right]$, a new operator. We now check the following:

**Claim:** $M_k^+ M_k |\Psi> = |\Psi>$.

**Proof:** Note that $<\Psi|\Psi> = 1$. Consider the case $k = 1$ as follows:

We have $M_1 = \left[2|\Psi><\Psi| - I\right]$, therefore,

$$M_1^+ M_1 |\Psi> = \left[2|\Psi><\Psi| - I\right]\left[2|\Psi><\Psi| - I\right]|\Psi>$$

$$= \left[4|\Psi><\Psi|\Psi><\Psi| - 2|\Psi><\Psi| - 2|\Psi><\Psi| + I\right]|\Psi>$$

$$= \left[4|\Psi><\Psi|\Psi><\Psi| - 4|\Psi><\Psi| + I\right]|\Psi>$$

$$= \left[4|\Psi><\Psi| - 4|\Psi><\Psi| + I\right]|\Psi>$$

$= |\Psi>$. Let us now consider the case $k = 2$ as follows:

$$M_2^+ M_2 |\Psi> = \left[4|\Psi><\Psi| - 3I\right]\left[4|\Psi><\Psi| - 3I\right]|\Psi>$$

$$= \left[16|\Psi><\Psi|\Psi><\Psi| - 12|\Psi><\Psi| - 12|\Psi><\Psi| + 9I\right]|\Psi>$$

$$= \left[16|\Psi><\Psi|\Psi><\Psi| - 24|\Psi><\Psi| + 9I\right]|\Psi>$$

$$= \left[16|\Psi><\Psi| - 24|\Psi><\Psi| + 9I\right]|\Psi>$$

$$= \left[16|\Psi> - 24|\Psi> + 9|\Psi>\right]|$$

$= |\Psi>$. On similar lines the general case also follows:

$$M_k^+ M_k |\Psi> = \left[2^k |\Psi><\Psi| - \left(2^k - 1\right)I\right]\left[2^k |\Psi><\Psi| - \left(2^k - 1\right)I\right]|\Psi>$$

$$= \left[2^{2k} |\Psi><\Psi| - 2(2^k(2^k - 1))|\Psi><\Psi| + (2^k - 1)(2^k - 1)I\right]|\Psi>$$

$$= \left[2^{2k} - 2^{2k} - 2^{2k} + 2^k + 2^k + 2^{2k} - 2^k - 2^k + 1\right]|\Psi>$$

$$= \left[2^{2k} - 2(2^{2k}) + 2(2^k) + 2^{2k} - 2(2^k) + 1\right]|\Psi>$$

$= |\Psi>$.

We now define the operator called the “total operator”, $T_k = M_k O$. We are now ready to discuss our second algorithm which requires **only a single step** to find the target state! Before we discuss the algorithm we state one important result which is used in this algorithm.

**Claim:** Let the initial wave function, $|\Psi>$, representing the quantum bag of indices be an equally weighted superposition of computational basis states, $|x>$, of length $k$ , i.e. the quantum bag contains $2^k$ indices. Also, let there be only one target state, $|t>$, then the target state, $|t>$, can be found, or reached, or attained, or achieved by just operating only once the operator $T_{(k-1)} = M_{(k-1)}O$ on this wave function, $|\Psi>$.

**Proof:** It is clear to see that $|\Psi>=\frac{1}{\sqrt{2^k}}\sum_{x=0}^{(2^k-1)}|x>$, where $|x>=|i_1 i_2 \cdots i_k>$

We now operate the operator $T_{(k-1)} = M_{(k-1)}O$ on the wave function $|\Psi>$. Thus we have

$$T_{(k-1)}|\Psi>=M_{(k-1)}O|\Psi>=M_{(k-1)}\left[|\Psi>-\frac{2}{\sqrt{2^k}}|t>\right]$$

$$=\left[(2^{(k-1)}|\Psi><\Psi|-(2^{(k-1)}-1)I)\right]\left[|\Psi>-\frac{2}{\sqrt{2^k}}|t>\right]$$

$$=2^{(k-1)}|\Psi>-2^{(k-1)}|\Psi>+|\Psi>-\frac{2^k}{\sqrt{2^k}}|\Psi><\Psi|t>+\frac{2(2^{(k-1)}-1)}{\sqrt{2^k}}|t>$$

$$=2^{(k-1)}|\Psi>-2^{(k-1)}|\Psi>+|\Psi>-\left(\frac{2^k}{\sqrt{2^k}}\right).\left(\frac{1}{\sqrt{2^k}}\right)|\Psi>+\frac{2(2^{(k-1)}}{\sqrt{2^k}}|t>$$

$$=2^{(k-1)}|\Psi>-2^{(k-1)}|\Psi>+|\Psi>-|\Psi>+\frac{2(2^{(k-1)}}{\sqrt{2^k}}|t>$$

$$=\frac{2(2^{(k-1)}-1)}{\sqrt{2^k}}|t>$$

Thus, we have got only the target state, $|t>$, with nonzero amplitude and all other basis states in $|\Psi>$ vanish, i.e. their amplitude becomes zero! Note that the amplitude of target state becomes large (in fact bigger than unity). This implies that the total "operator" is not unitary, since the action of unitary operator on a vector preserves its length and here the chosen vector (wave function $|\Psi>$) is of unit length ( $\||\Psi>\|^2=1$ ).

We now proceed to formally discuss the steps of the algorithm which consists of just applying the appropriate "total operator" on the wave function, $|\Psi>$, representing the given quantum bag of indices containing a single target index.
Thus, let the given bag of indices contains $N=2^n$ elements. We will prepare the quantum bag in terms of the wave function, $|\Psi>$, as follows in the following

## Steps of the algorithm:

**(i)** We consider a quantum state containing $n$ qubits, all initialized to zero, i.e. the state

$$|00\cdots0>=|0>^{\otimes n}.$$

**(ii)** We apply Hadamard transform to all the $n$ qubits to get

$$|\Psi>=H^{\otimes n}|0>^{\otimes n}=\frac{1}{\sqrt{2^n}}\sum_{x=0}^{2^n-1}|x>$$ . Clearly, $|t>\in|\Psi>$.

**(iii)** Since the size of the data is $N=2^n$ so we choose $T_{(n-1)}=M_{(n-1)}O$ as our "total operator" to operate on the wave function $|\Psi>$.

**(iv)** We carry out the action of the chosen operator on the wave function only in terms of the target state itself!

$$T_{(n-1)}|\Psi>=M_{(n-1)}O|\Psi>=M_{(n-1)}\left[|\Psi>-\frac{2}{\sqrt{2^n}}|t>\right]$$

$$=\left[(2^{(n-1)}|\Psi><\Psi|-(2^{(n-1)}-1)I)\right]\left[|\Psi>-\frac{2}{\sqrt{2^n}}|t>\right]$$

$$=\frac{2(2^{(n-1)}-1)}{\sqrt{2^n}}|t>.$$

□

Thus, it is clear that if we carry out measurement then we will get the target state, $|t>$, with probability one!! Thus, this algorithm assures us to obtain the target state with 100% guarantee!!!

### 3.6 The Third New Quantum Search Algorithm implying *P* = *NP* :

Again, Let $x$ , $0\le x\le(N-1)$, be an element in the unstructured database of size $N$ . Let $N=2^n$ hence $0\le x\le(2^n-1)$. Our aim in the unstructured database search problem is to locate and pick out the target index, $t\in[0,(2^n-1)]$. Also, Suppose we have a 1-YES quantum oracle defined in terms of operator, $O$, which performs the operation $O|x>=(-1)^{f_t(x)}|x>$, where as mentioned previously $f_t(x)$ is defined such that $f_t(x)=1$ if $x=t$ and $f_t(x)=0$ if $x\ne t$ . where $f_t(x)=1$ stands for YES and $f_t(x)=0$ stands for NO. As seen previously, the operator $O$ is unitary. Thus, everything is same as it was in previous algorithms. In this new quantum algorithm we will be doing non-unitary quantum computation, i.e. the operator we will be using to achieve the task of enhancing the amplitude of the target state, $|t>$, as is done in the previous two algorithms is non-unitary. This algorithm also works in just a single step, i.e. it enhances the amplitude of the target state to its full in just one operation of the non-unitary operator chosen for this algorithm.

### Steps of the algorithm:

**(i)** We consider a quantum state containing $n$ qubits, all initialized to zero, i.e. the state $|00\cdots 0>=|0>^{\otimes n}$.

**(ii)** We apply Hadamard transform to all the $n$ qubits to get

$$|\Psi>=H^{\otimes n}|0>^{\otimes n}=\frac{1}{\sqrt{2^n}}\sum_{x=0}^{2^n-1}|x>$$. Clearly, $|t>\in|\Psi>$.

**(iii)** We apply non-unitary operator, $A=\frac{\sqrt{N}}{2}[I-O]$, on the wave function $|\Psi>$. We get

$$A|\Psi>=\frac{\sqrt{N}}{2}[I-O]|\Psi>=\frac{\sqrt{N}}{2}\left[|\Psi>-|\Psi>+\frac{2}{\sqrt{N}}|t>\right]=|t>$$

.

□

Thus, we have seen that by the action of non-unitary operator, $A=\frac{\sqrt{N}}{2}[I-O]$, and carry out the measurement then we will get the target state, $|t>$, with probability one! The thing to be seen is whether it is possible to build quantum circuit which will perform the action of the non-unitary operator, $A=\frac{\sqrt{N}}{2}[I-O]$.

## 4 Remarks:

**Remark 1:** It is clear to see that as the algorithm proceeds we get at each iteration the bag containing proper range of indices to which target index belongs, i.e. we get during each of the iterations a proper quantum bag reduced to half in size, in terms of a quantum state which has nonzero overlap with target state. Thus, as we proceed at an intermediate stage we reach at a wave function,

$$| \Psi_{i_1 i_2 \cdots i_k} >= | i_1 i_2 \cdots i_k > \otimes H^{\otimes(n-k)} | 00 \cdots 0 >$$

which has nonzero overlap with the target state, $| t >$. We then divide the quantum bag into two new quantum bags, i.e. construct two new states out of which only one will have nonzero overlap with the target state, $| t >$, to be determined by taking inner product with any one of these two newly prepared quantum states. Thus, the new quantum states constructed from consideration of the earlier reached above mentioned quantum state will be

$$| \Psi_{i_1 i_2 \cdots i_k 0} >= | i_1 i_2 \cdots i_k > \otimes | 0 > \otimes H^{\otimes(n-k-1)} | 00 \cdots 0 >$$

and,

$$| \Psi_{i_1 i_2 \cdots i_k 1} >= | i_1 i_2 \cdots i_k > \otimes | 1 > \otimes H^{\otimes(n-k-1)} | 00 \cdots 0 >.$$

**Remark 2:** It is interesting to see that the amplitude of each state in the equally weighted superposition of states (including target state) is initially equal to $\frac{1}{\sqrt{N}}$. This state represents the initial quantum bag. After first iteration of the size of the quantum bag reduces to half and this size reaches finally to unity, i.e. finally (at the $n^{th}$ iteration) the quantum bag will contain only the target state itself! Therefore, after first iteration the amplitude of each state in the equally weighted superposition of states becomes $\sqrt{\frac{2}{N}}$. The amplitude of each state including target state in the superposition changes in the successive iterations as follows:

$$\frac{1}{\sqrt{N}} \rightarrow \sqrt{\frac{2}{N}} \rightarrow \sqrt{\frac{2^2}{N}} \rightarrow \cdots \rightarrow \sqrt{\frac{2^j}{N}} \rightarrow \cdots \rightarrow 1$$

**Remark 3:** It is clear to see that in $n = LogN$ iterations we will attain the target state, i.e. in the final quantum bag, after carrying out $n = LogN$ iterations, will contain only the target state $| t >$ itself which will lead to the value of inner product equal to unity.

**Remark 4:** It is important to note that actually in each iteration of the algorithm we are getting one bit of the target state. i.e. if the target state is $| t >= | t_1 t_2 \cdots t_j \cdots t_n >= | t_1 > \otimes | t_2 > \otimes \cdots \otimes | t_j > \otimes \cdots \otimes | t_n >$ then in first iteration we determine the first bit namely, $| t_1 >$, in the successive iterations we determine $| t_2 >, \cdots, | t_j >, \cdots, | t_n >$. Thus in $n = LogN$ iterations we will be able to determine the target state, $| t >$, completely.

**Remark 5:** Alternatively, instead of one oracle we may define implicitly *n* number of oracles, $f_t^i(x)$, which gives rise to *n* number of target states $|t^i>=|t_i t_{(i+1)}\cdots t_n>$. Clearly, $f_t^1(x)=f_t(x)$ and it gives rise to target state $|t^1>=|t>$. Further, by finding the nonzero inner product between the inner products taken that of $|t^i>$ with any one of the wave functions, $|\Psi_{i0}>$ and $|\Psi_{i1}>$ that we build, namely, $|\Psi_{i0}>=|0>\otimes H^{(n-i-1)}|00\cdots0>$ and $|\Psi_{i1}>=|1>\otimes H^{(n-i-1)}|00\cdots0>$, we can determine separately each bit $|t_i>$ of the target state $|t>$ and then build it as $|t>=|t_1t_2\cdots t_j\cdots t_n>$.

**Remark 6:** As far as the value of inner product is concerned we are only interested to know whether it is zero or nonzero, and we are not interested in its exact value. Therefore we can use the existing quantum algorithm [2] to evaluate the inner product with complexity $\sim O(LogN)$. Since our new quantum algorithm requires $LogN$ steps to reach the desired target state and each iterative step requires to find out one inner product which again takes time $\sim O(LogN)$ therefore, our new quantum search algorithm is of the order $\sim O((LogN)^2)$.

**Remark 7:** For a typical ***NP-Complete*** **problem** in which one has to find an assignment of one of the $b$ values to each of the $C$ variables, the number of candidate solutions, $N=b^C$, grows exponentially with $C$. Hence, the classical algorithm for unstructured search would therefore take time of the order, $\sim O(b^C)$, to find desired solution (as the target state) e.g. minimum weight Hamiltonian circuit among the all possible Hamiltonian circuits as a solution for the traveling salesman problem, whereas the Grover's quantum algorithm [9] would take a time of the order, $\sim O(b^{\frac{C}{2}})$. But from the complexity of the order $\sim O((LogN)^2)$ that we get for our quantum search algorithm it is easy to check that our quantum search algorithm will takes time of the order, $\sim O(b^{(LogC)^2})$, thus an impressive (exponential) speedup over existing classical or quantum algorithm. We thus have managed $P=NP$ using our new quantum search algorithm.

**Example 1:** Let the bag of indices contains numbers $\{0,1,2,\cdots,15\}$ and let the target element, $t=11$. We begin with the wave function, $|\Psi>$, namely, $|\Psi>=H^{\otimes 4}|0000>$ which contains the target state, $|t>=|11>=|1011>$. We now follow the steps of the algorithm:

Cleary, $<t|\Psi>\neq 0$, therefore, we divide quantum bag represented by $|\Psi>$ into two bags, represented by $|\Psi_0>$ and $|\Psi_1>$, where

$|\Psi_0>=|0>\otimes H^{\otimes 3}|000>$, and $|\Psi_1>=|1>\otimes H^{\otimes 3}|000>$.

Cleary, $<t|\Psi_1>\neq 0$, therefore, we further divide quantum bag represented by $|\Psi_1>$ into two bags, represented by $|\Psi_{10}>$ and $|\Psi_{11}>$, where

$$|\Psi_{10}>=|1>\otimes|0>\otimes H^{\otimes 2}|00>\text{ , and}$$

$$|\Psi_{11}>=|1>\otimes|1>\otimes H^{\otimes 2}|00>.$$

Cleary, $<t|\Psi_{10}>\neq 0$, therefore, we further divide quantum bag represented by $|\Psi_{10}>$ into two bags, represented by $|\Psi_{100}>$ and $|\Psi_{101}>$, where

$$|\Psi_{100}>=|1>\otimes|0>\otimes|0>\otimes H|0>\text{, and}$$

$$|\Psi_{101}>=|1>\otimes|0>\otimes|1>\otimes H|0>.$$

Cleary, $<t|\Psi_{101}>\neq 0$, therefore, we further divide quantum bag represented by $|\Psi_{101}>$ into two bags, represented by $|\Psi_{1010}>$ and $|\Psi_{1011}>$, where

$$|\Psi_{1010}>=|1>\otimes|0>\otimes|1>\otimes|0>\text{, and}$$

$$|\Psi_{1011}>=|1>\otimes|0>\otimes|1>\otimes|1>.$$

Clearly, $<t|\Psi_{1011}>\neq 0$, and in fact $<t|\Psi_{1011}>=1$, therefore,
We have located (reached to) the desired target state, $|t>=|11>=|1011>$, present in the given database (initial quantum bag containing target) in terms of the superposition state, $|\Psi>=H^{\otimes 4}|0000>$.

**Example 2:** Let the bag of indices contains numbers $\{0,1,2,\cdots,7\}$ and let the target element, $t=3$. We begin with the wave function, $|\Psi>$, namely, $|\Psi>=H^{\otimes 3}|000>$ which contains the target state, $|t>=|3>=|011>$.
Carrying out **step (iv) of the second algorithm** we have

$$\left[(2^2|\Psi><\Psi|-(2^2-1)I)\right]\left[O\right]\left[|\Psi>\right]$$

$$=\left[(2^2|\Psi><\Psi|-(2^2-1)I)\right]\left[|\Psi>-\frac{2}{\sqrt{2^3}}|t>\right]$$

$$=\left[4|\Psi><\Psi|-3I\right]\left[|\Psi>-\frac{1}{\sqrt{2}}|t>\right]$$

$$=4|\Psi>-3|\Psi>-\frac{4}{\sqrt{2}}|\Psi><\Psi|t>+\frac{3}{\sqrt{2}}|t>$$

$$= 4\,|\,\Psi> - 3\,|\,\Psi> - \frac{4}{\sqrt{2}}\,|\,\Psi>.\frac{1}{2\sqrt{2}} + \frac{3}{\sqrt{2}}\,|\,t> = \frac{3}{\sqrt{2}}\,|\,t>.$$

**Example 3:** Consider example same as Example 1 above. We solve it now using second algorithm: Let the bag of indices contains numbers $\{0,1,2,\cdots,15\}$ and let the target element, $t = 11$. We begin with the wave function, $|\,\Psi>$, namely, $|\,\Psi> = H^{\otimes 4}\,|\,0000>$ which contains the target state, $|\,t> = |\,11> = |\,1011>$. Carrying out **step (iv) of the second algorithm** we have

$$\left[(2^3\,|\,\Psi><\Psi\,| - (2^3-1)I)\right]\left[O\right]\left[|\,\Psi>\right]$$

$$= \left[(2^3\,|\,\Psi><\Psi\,| - (2^3-1)I)\right]\left[|\,\Psi> - \frac{2}{\sqrt{2^4}}\,|\,t>\right]$$

$$= \left[8\,|\,\Psi><\Psi\,| - 7I\right]\left[|\,\Psi> - \frac{1}{2}\,|\,t>\right]$$

$$= 8\,|\,\Psi> - 7\,|\,\Psi> - 4\,|\,\Psi><\Psi\,|\,t> + \frac{7}{2}\,|\,t>$$

$$= 8\,|\,\Psi> - 7\,|\,\Psi> - 4\,|\,\Psi>.\frac{1}{4} + \frac{7}{2}\,|\,t> = \frac{7}{2}\,|\,t>.$$

**Example 4:** Consider same example above. We now solve it using third algorithm: Let the bag of indices contains numbers $\{0,1,2,\cdots,15\}$ and let the target element, $t = 11$. We begin with the wave function, $|\,\Psi>$, namely, $|\,\Psi> = H^{\otimes 4}\,|\,0000>$ which contains the target state, $|\,t> = |\,11> = |\,1011>$. We carry out **step (iii) of the third algorithm, i.e.** we apply non-unitary operator, $A = \frac{\sqrt{N}}{2}\left[I - O\right]$, on the wave function $|\,\Psi>$. This gives $A = \frac{\sqrt{N}}{2}\left[I - O\right] = \frac{4}{2}\left[I - O\right]H^{\otimes 4}\,|\,0000> = |\,t>$.

## 5. A novel efficient quantum algorithm for solving K-SAT problem, K ≥3, on a quantum computer:

In this section we proceed to develop a novel efficient quantum algorithm for solving K-SAT problem, K ≥3, on a quantum computer which further can be led to develop not only an efficient quantum algorithm for solving traveling salesman problem via spending an additional polynomial cost of translation [12].This efficient quantum algorithm for solving K-SAT problem, K $\geq$3, on a quantum computer not only leads to an efficient quantum algorithm for solving traveling salesman problem but in fact it leads to getting efficient quantum algorithms for solving all NP-complete problems via spending additional polynomial cost of translation as mentioned [12].

In computer science, a decision problem is a problem with a yes or no answer. A decision problem is in NP (which stands for Non-deterministic Polynomial time) if a yes answer can be verified efficiently, i.e. in a time that grows no faster than a polynomial in the size of the problem. Hence, loosely speaking, the problems in NP are those such that if you happened to guess the solution correctly (this is the non-deterministic aspect) then you could verify the solution efficiently (this is the polynomial aspect). Hence the name Non-deterministic Polynomial time. A decision problem is NP-complete if it lies in the complexity class NP and all other problems in NP can be reduced to it. Thus the NP-complete problems are the only problems we need to study to understand the computational resources needed to solve all of the problems in NP [10]. Hence, the NP-complete problems have a special place in complexity theory. In addition to their ubiquity, NP-complete problems share a fortuitous kinship: any NP-complete problem can be mapped into any other NP-complete problem using only polynomial resources. Thus, if an efficient algorithm were found that can solve one type of NP-complete problem efficiently, this would immediately lead to efficient algorithms for solving all NP-complete problems (up to the polynomial cost of translation). So, in some sense, the three thousand or so different NP-complete problems are really the same problem in disguise. It is therefore sufficient to focus on any one NP-complete problem, for any progress made in solving that problem is likely to be applicable to all the other NP-complete problems too, so long as you do not mind paying the polynomial cost of translation.

Although a solution to an NP-complete problem can be *verified* "quickly", there is no known way to *find* a solution quickly. That is, the time required to solve the problem using any currently known algorithm increases rapidly as the size of the problem grows.

As a consequence, determining whether it is possible to solve these problems quickly, called the P versus NP problem, is one of the fundamental unsolved problems in computer science today.

P versus NP problem corresponds to the question of whether all so-called NP problems are actually P problems. A P problem is one that can be solved in "polynomial time," which means that an algorithm exists for its solution such that the number of steps in the algorithm is bounded by a polynomial function of $n$, where $n$ corresponds to the length of the input for the problem. Thus, P problems are said to be easy, or tractable. While a method for computing the solutions to NP-complete problems quickly remains undiscovered, computer scientists and programmers still frequently encounter NP-complete problems. NP-complete problems are often addressed by using heuristic methods and approximation algorithms. Thus, NP-complete problems form a class of computational problems for which no efficient solution algorithm has been found so far. Many significant computer-science problems belong to this class, e.g., the satisfiability problems (K-SAT), traveling salesman problem, and graph-covering problems.

It is not known whether any polynomial-time (efficient) algorithms will ever be found for NP-complete problems, and determining whether these problems are tractable or intractable remains one of the most important questions in theoretical computer science. When an NP-complete problem must be solved, one approach is to use a polynomial algorithm to approximate the solution; the answer thus obtained will not necessarily be optimal but will be reasonably close.

A famous example of an NP-complete problem is the traveling salesman problem, which has wide applications in the optimization of transportation schedules. Determining whether NP-complete problems are tractable or intractable remains one of the most important questions in theoretical computer science. Such a discovery would prove that P = NP = NP-complete and revolutionize many fields in computer science and mathematics.

Proving the result P = NP will lead to significant theoretical and practical consequences. This result will make huge impact on various branches of science and technology.

To develop a new efficient quantum algorithm in order to solve some NP-Complete problems we choose the famous K-SAT problem, $K \geq 3$, and give an efficient quantum algorithm to solve this problem efficiently. It is important to note that actually this so called K-SAT problem was the first problem shown to be NP-Complete [1]. The chosen K-SAT problem, for which $K \geq 3$, will be involving in all n variables, defined by the Boolean function, $\phi$, which is having a unique assignment, $x = t$, for which $\phi(t) = 1$ and $\phi(x) = 0$ for all other assignments, $x( \neq t)$. NP-Complete problems share a fortuitous kinship: any NP-Complete problem can be mapped into any other NP-complete problem using only polynomial resources [12]. Thus, if an efficient algorithm were found that can solve one type of NP-Complete problem efficiently, this would immediately lead to efficient algorithms for all the other NP-complete problems (up to the polynomial cost of translation).
NP-complete problems form a class of computational problems for which no efficient solution algorithm has been found so far. Many significant computer-science problems belong to this class, e.g., the satisfiability problems (K-SAT), traveling salesman problem, and graph-covering problems. Proving the result P = NP will lead to significant theoretical and practical consequences. This result will make huge impact on various branches of science and technology. A quantum computer is a computer that exploits quantum mechanical phenomena. A scalable quantum computer could perform some calculations exponentially faster than any modern "classical" computer. It has been shown that some quantum algorithms, e.g. Shor's quantum algorithm for factoring a number into its prime factors, are exponentially more efficient than the best-known classical algorithms for doing the same task. A large-scale quantum computer could in theory solve computational problems unsolvable by a classical computer in any reasonable amount of time. This concept of extra ability has been called "quantum supremacy". We have discovered a quantum algorithm for solving K-SAT problem, $K \geq 3$, which is an NP-Complete problem. This newly discovered quantum algorithm for solving K-SAT problem, $K \geq 3$, is exponentially more efficient than the best-known classical algorithms for doing the same task. An NP-complete problem can be mapped into any other NP-complete problem using only polynomial resources. Using this well-known result the efficient quantum algorithms will be obtained for efficiently solving all known (three thousand plus) NP-complete problems by mapping them on the K-SAT problem, $K \geq 3$, using only polynomial cost of translation. Thus, this patent offers an efficient quantum algorithm for efficiently solving K-SAT problem, $K \geq 3$, which is the first problem shown to be NP-Complete, and by implication every other problem among the NP-complete problems by bearing the polynomial cost of translation.

## 5.1 A Quantum Algorithm for Solving K-SAT Problem:

We choose K-SAT problem, $K \geq 3$, among the NP-complete problems to solve efficiently, which is involving n variables, and is defined by the Boolean function, $\phi$. This K-SAT problem, $K \geq 3$, was the first problem shown to be NP-complete. We choose the K-SAT problem of most difficult type, i.e. which has a unique correct assignment, $x = t$, that makes it true, i.e. $\phi: \{0, 1\}^n \rightarrow \{0, 1\}$ such that for $t = (t_1, t_2, .... , t_n)$, $t_i \varepsilon \{0, 1\}$, $\phi(t) = 1$, and $\phi(x) = 0$ for all $x \neq t$, $x = (x_1, x_2, .... , x_n)$, $x_i \varepsilon \{0, 1\}$. We propose a novel polynomial time quantum algorithm for finding this correct assignment $t = (t_1, t_2, .... , t_{(n-1)}, t_n)$ making the Boolean function true in terms of corresponding quantum state $|t> = |t_1 t_2 t_3 .... t_{(n-1)} t_n>$ out of all possible $2^n$ assignments $x = (x_1, x_2, .... , x_n)$, $x_i \varepsilon \{0, 1\}$. The proposed quantum algorithm determines $|t> = |t_1 t_2 t_3 t_4 .... t_{(n-1)} t_n>$ by performing n iterative steps. The algorithm determines the single qubit, $|t_i>$, in the i-th iterative step, $i \varepsilon \{1, 2, .... , n\}$, which makes a part of the entire solution state, $|t>$, containing n qubits for K-SAT problem under consideration. Thus, in the first iterative step the algorithm determines the first qubit, $|t_1>$, and it then determines, in the successive iterations, the other qubits, $|t_2>$, $|t_3>$, .... , $|t_{(n-1)}>$, $|t_n>$ and thus it completely determines the solution state, $|t>$.

The algorithm begins with preparing two quantum registers $|\psi_0>$ and $|\psi_1>$ such that $|\psi_0>$ contains superposition of the computational basis states corresponding to first fifty percent configurations for the Boolean function, $\phi$, so that obviously the first qubit of each computational basis state belonging to this register representing an assignment will be $|0>$, The other register $|\psi_1>$ contains superposition of the remaining computational basis states corresponding to the next fifty percent configurations for the Boolean function, $\phi$, so that obviously the first qubit of each computational basis state belonging to this register representing an assignment will be $|1>$. Clearly, the target state, $|t>$, will be present either in $|\psi_0>$ or in $|\psi_1>$. We then operate on these registers, the unitary operator, $U_\phi$, corresponding to the Boolean function $\phi$. The action of $U_\phi$, on $|\psi_0>$ will invert the phase of the solution state $|t>$ if and only if it is present in this register and if the solution state $|t>$ is not present in the register $|\psi_0>$ then in this case $|\psi_0>$ will remain unaffected (invariant) under the action of $U_\phi$. Same is true for the register $|\psi_1>$.

Case 1: If $|\psi_0>$ remains invariant under the action of $U_\phi$, i.e. $U_\phi |\psi_0> = |\psi_0>$ then this implies that $|t>$ will then belong to register $|\psi_1>$ and clearly in this case $|t_1> = |1>$ since all the states present in the superposition of states belonging to this register $|\psi_1>$ begin with the single qubit state $|1>$.

How to check whether $U_\phi |\psi_0> = |\psi_0>$ or not?

We proceed as follows: In order to create register, $|\psi_0>$, one first creates the state, $|0>^{(\otimes n)}$, where $|0>^{(\otimes n)}$ represents the tensor product of n copies of the state, $|0>$, and then operates on this state the unitary operator, $I \otimes H^{(\otimes(n-1)}$ where I denotes the identity operator and $H^{(\otimes(n-1)}$ denotes the tensor product of (n-1) copies of Hadamard operator. Therefore, if one will apply once more the same (unitary) operator, $I \otimes H^{(\otimes(n-1)}$, on thus formed register, $|\psi_0>$, then one will get back the original state, $|0>^{(\otimes n)}$, where one has started. Now, when Case (i) is true, i.e. when $|\psi_0>$ remains invariant under the action of the operator, $U\phi$, i.e. when $U_\phi |\psi_0> = |\psi_0>$ then it is clear to see that one will get back the original starting state, $|0>^{(\otimes n)}$, when one will operate the same above mentioned operator, $I \otimes H^{(\otimes(n-1)}$, on $U_\phi |\psi_0>$ and this will verify that Case (i) is true.

Case 2: Instead, if register $|\psi_1>$ remains invariant under the action of unitary operator $U_\phi$, i.e. if $U\phi|\psi_1>. = |\psi_1>.$ then it implies that $|t>$ belongs to register $|\psi_0>$. Further, in this case $|t_1> = |0>$ since all the states present in the superposition of states belonging to this

register $|\psi_0>$ begin with the single qubit state $|0>$.

How to check whether $U_\phi |\psi_1> = |\psi_1>$ or not?

In order to create register, $|\psi_1>$, one first creates the state, $|0>^{(\otimes n)}$, where $|0>^{(\otimes n)}$ represents the tensor product of n copies of the state, $|0>$, and then operates on this state the unitary operator, $X \otimes H^{(\otimes(n-1)}$ where I denotes the NOT operator and $H^{(\otimes(n-1)}$ denotes the tensor product of (n-1) copies of Hadamard operator. Therefore, if one will apply once more the same (unitary) operator, $X \otimes H^{(\otimes(n-1)}$, on thus formed register, $|\psi_1>$, then one will get back the original state, $|0>^{(\otimes n)}$, where one has started. Now, when Case (ii) is true, i.e. when $|\psi_1>$ remains invariant under the action of the operator, $U\phi$, i.e. when $U_\phi |\psi_1> = |\psi_1>$ then it is clear to see that one will get back the original starting state, $|0>^{(\otimes n)}$, when one will operate the same above mentioned operator, $X \otimes H^{(\otimes(n-1)}$, on $U_\phi |\psi_1>$ and this will verify that Case (ii) is true.

Making use of the above cases we determine $|t_1>$ in the first iteration as described in the next section.

(i) When we get $|t_1> = |0>$ it implies that $|t>$ belongs to register $|\psi_0>$ and in this case we divide the register $|\psi_0>$ into two equal parts. This is achieved essentially by creating two new registers $|\psi_{00}>$ and $|\psi_{01}>$ leading to the determination of $|t_2>$ as described in the next section.

(ii) On the other hand when we get $|t_1> = |1>$ it implies that $|t>$ belongs to register $|\psi_1>$ and we divide the register $|\psi_1>$ into two equal parts. This is achieved essentially by creating two new registers $|\psi_{10}>$ and $|\psi_{11}>$ leading to the determination of $|t_2>$ as described in the next section.

We will see in the next section that by repeating the same procedure of dividing the appropriate registers through creating two new registers we determine other single qubits $|t_3>$, $|t_4>$, .... , $|t_n>$ etc. and thus using in total n such iterations we will determine the solution state, $|t>$.

Now, the effect of the action of $U_\phi$ on the register containing superposition of entire computational basis is as follows:

Let $|\psi>$ be a quantum register made up of equally weighted superposition of all the computational basis states such that each one of this computational basis state represents a possible assignment for the given K-SAT problem. Thus,

$|\psi> = H^{(\otimes n)}|0>^{(\otimes n)}$

where $H^{(\otimes n)}$ stands for tensor product of n copies of Hadamard gate and $|0>^{(\otimes n)}$ stands for tensor product of n copies of state $|0>$. Let $U_\phi$ denotes the unitary operator corresponding to the Boolean function, $\phi$.

Let $|\Omega> = H|1> = 1/\sqrt{2}(|0> - |1>)$

The action of $U_\phi$ is as follows:

$U_\phi : |x>|\Omega> \rightarrow |x>|\Omega + \phi(x) >$

where $\phi(x) = 1$, if $x = t$, the only assignment that makes the Boolean function $\phi$ true, and $\phi(x) = 0$ otherwise. $| \Omega >$ is the ancilla qubit, and $|\Omega + \phi(x)>$ is the exclusive-OR of the bit value of the ancilla and the bit value that is output from our black-box function $\phi(x)$.

Thus, when the unitary operator $U_\phi$ will be operated on the quantum register $|\psi>$ containing superposition of all the computational basis states, representing together all possible assignments for the given K-SAT problem, the phase associated with the only computational basis state, $|x> = |t>$, which is the only solution for the K-SAT problem under consideration, i.e. for which $\phi(t) = 1$, will be inverted. In short, for assignment $x = t$ for which the Boolean function $\phi$ is true, giving rise to $\phi(t) = 1$, we will have $U_\phi(|t>) = -|t>$, and for all other assignments $x(\neq t)$ we will have $U_\phi(|x>) = |x>$.

Thus, when we compute the given Boolean function, $\phi$, on this register, $|\psi>$, i.e. when we operate the unitary operator, $U_\phi$, corresponding to the given Boolean function, $\phi$, on this register $|\psi>$ containing superposition of all the computational basis states, representing together all possible assignments for the given K-SAT problem and evaluate $U_\phi|\psi>$ it will cause the inversion of the phase of the only computational basis state, $|t> = |t_1 t_2 \ldots t_{(n-1)} t_n>$, which corresponds to the unique solution $t = (t_1, t_2, \ldots, t_n)$ to the given K-SAT problem at hand.

## 5.2. Steps of the New Efficient Quantum Algorithm for Solving K-SAT Problem, K ≥ 3:

1. We prepare the following two quantum registers, $|\Psi_0>$ and $|\Psi_1>$ such that the first register contains the superposition of first fifty percent of the computational basis states and the second register contains the remaining next fifty percent of the computational basis states, i.e. we construct two quantum registers $|\Psi_0>$ and $|\Psi_1>$ as follows:
   $|\Psi_0> = |0> \otimes H^{(\otimes(n-1))}|0^{(\otimes(n-1))}>$
   and,
   $|\Psi_1> = |1> \otimes H^{(\otimes(n-1))}|0^{(\otimes(n-1))}>$.
   Clearly, the computational basis state $|t>$ that corresponds to the unique solution, $t = (t_1, t_2, \cdot\cdot\cdot, t_n)$ to the given K−SAT problem will belong to some one of the registers $|\Psi_0>$ and $|\Psi_1>$, i.e. either $|t> \in |\Psi_0>$ or $|t> \in |\Psi_1>$.

2. We now apply the unitary operator, $U_\Phi$, that computes $\Phi$ on these registers. It is important to note that if $U_\Phi|\Psi_0> = |\Psi_0>$, i.e. if $||\Psi_0> - U_\Phi|\Psi_0>| = 0$, or equivalently, if $|<\Psi_0|U_\Phi|\Psi_0>| = 1$ then $U_\Phi|\Psi_0>$ will be a "separable" state and in this case $|t> \in |\Psi_1>$ and therefore in this case we will get $U_\Phi|\Psi_1> = |\Psi_1> - 2/\sqrt{2}^{(n-1)}|t>$, an "entangled" state and as a consequence we will have $||\Psi_1> - U_\Phi|\Psi_1>| = |2/\sqrt{2}^{(n-1)}|t>|$, or equivalently, $|<\Psi_1|U_\Phi|\Psi_1>| < 1$. On the other hand if $U_\Phi|\Psi_1> = |\Psi_1>$, i.e. if $||\Psi_1> - U_\Phi|\Psi_1 i| = 0$, or equivalently, if $|<\Psi_1|U_\Phi|\Psi_1>| = 1$ then $U_\Phi|\Psi_1>$ will be a "separable" state and in this case $|t> \in |\Psi_0>$ and therefore in this case we will get $U_\Phi|\Psi_0> = |\Psi_0> - 2/\sqrt{2}\,(n-1)\,|t>$, an "entangled" state and as a consequence we will have $||\Psi_0> - U_\Phi|\Psi_0>| = |2/\sqrt{2}\,(n-1)\,|t>|$, or equivalently, $|<\Psi_0|U_\Phi|\Psi_0>| < 1$. Thus, the effect the action of $U_\Phi$ on the states $|\Psi_0>$ and $|\Psi_1>$ is totally different. Rather it is diametrically opposite. Between the states $|\Psi_0>$ and $|\Psi_1>$ the state to which $|t>$ does not belong remains invariant under the action of $U_\Phi$ and it remains as it was initially i.e. a fully "separable" state made up of the product of linear (single qubit) factors while the other state to which $|t>$ belongs becomes a fully "entangled" state not having any linear (single qubit) factor. Using this marked difference between the resulting states under the action of the operator $U_\Phi$ on the states $|\Psi_0>$ and $|\Psi_1>$ i.e. using this marked difference between the states $U_\Phi|\Psi_0>$ and $U_\Phi|\Psi_1>$ we can determine the register to which $|t>$ belongs (and the register to which $|t>$ belongs is to be used for the next iteration to find $|t2>$ as is described below).

3. It is now clear to see that after the action of the unitary operator, $U_\Phi$, one and only one of the quantum registers, either $|\Psi_0>$ or $|\Psi_1>$ given above in 1, remains invariant and so a "separable" state while the other register (to which the solution state $|t>$ belongs) becomes fully "entangled". In fact after the action of the unitary operator, $U_\Phi$, on the two new quantum registers that arise at each iterative step, where one of the registers contains the solution state $|t>$ while the other register does not, form respectively a fully "entangled" state not having any linear (single qubit) factor and a fully "separable" state made up of the product of linear (single qubit) factors and so these are fully "distinguishable" quantum states. Thus:

   Case (i): Without loss of generality, if $U_\Phi|\Psi_0> = |\Psi_0>$ then $|t> \in |\Psi_1>$, and $|t_1> = |1>$.

   Case (ii): Without loss of generality, if $U_\Phi|\Psi_1> = |\Psi_1>$ then $|t> \in |\Psi_0>$, and $|t_1> = |0>$.

   It is important to note that in order to create register, $|\psi_0>$, one first creates the state, $|0>^{(\otimes(n))}$ and then operates on this state the operator, $I \otimes H^{(\otimes(n-1))}$where $I$ denotes the identity operator and $H$ denotes the Hadamard operator. Therefore, if one will apply once more the same operator, $I \otimes H^{(\otimes(n-1))}$ on thus formed register, $|\psi_0>$, then one will get back the original state, $|0>^{(\otimes(n))}$, where one has started. Now, when Case (i) is true, i.e. when $|\psi_0>$ remains invariant under the action of the operator, $U_\Phi$ , i.e. when $U_\Phi|\Psi_0> = |\Psi_0>$ then it is clear to see that one will get back the original starting state, $|0>^{(\otimes(n))}$, when one will operate the same above mentioned operator, $I \otimes H^{(\otimes(n-1))}$ on $U_\Phi|\Psi_0>$ and this will verify that Case (i) is true.

On similar lines one can note that in order to create register, $|\psi_1>$, one first creates the state, $|0>^{(\otimes(n))}$, and then operates on this state the operator, $X\otimes H^{(\otimes(n-1))}$, where $X$ denotes the NOT operator and $H$ denotes the Hadamard operator. Therefore, if one will apply once more the same operator, $X\otimes H^{(\otimes(n-1))}$, on thus formed register, $|0>^{(\otimes(n))}$, then one will get back the original state, $|\psi_1>$, where one has started. Now, when Case (ii) is true, i.e. when $|\psi_1>$ remains invariant under the action of the operator, $U_\Phi$ , i.e. when $U_\Phi|\Psi_1> = |\Psi_1>$ then it is clear to see that one will get back the original starting state, $|0>^{(\otimes(n))}$, when one will operate the same above mentioned operator, $X\otimes H^{(\otimes(n-1))}$, on $U_\Phi|\Psi_1>$ and this will verify that Case (ii) is true.

4. Suppose Case (ii) given in 3 above is true: Then we have $t \in [0,(2(n-1) - 1)]$, and in this case we divide the register $|\Psi_0>$ into two equal parts, i.e. we create two new registers as follows: We construct two new quantum registers $|\Psi_{00}>$ and $|\Psi_{01}>$ such that $|\Psi_{00}>$ will represent the register that contains computational basis states corresponding to all those indices, x, such that $0 \leq x \leq (2(n-2) - 1)$ and $|\Psi_{01}>$ will represent the register that contains computational basis states corresponding to all those indices, x, such that $2(n-2) \leq x \leq (2(n-1) - 1)$. Thus, we have
$|\Psi_{00}> = |0> \otimes |0> \otimes H^{(\otimes(n-2))}|0>^{(\otimes(n-2))}$
and,
$|\Psi_{01}> = |0> \otimes |1> \otimes H^{(\otimes(n-2))}|0>^{(\otimes(n-2))}$
Now, clearly: either $|t> \in |\Psi_{00}>$ or $|t> \in |\Psi_{01}>$, and this is to be determined in exactly similar way as described above in 3 by doing the action of the the unitary operator, $U_\Phi$, on these quantum registers, $|\Psi_{00}>$ and $|\Psi_{01}>$, which will lead to creation of fully "distinguishable" states, $U_\Phi|\Psi_{00}>$ and $U_\Phi|\Psi_{01}>$.
Thus, we have two cases again:
Case (i): Without loss of generality, if $U_\Phi|\Psi_{00}> = |\Psi_{00}>$ then $|t> \in |\Psi_{01}>$, and $|t_2> = |1>$.
Case (ii): Without loss of generality, if $U_\Phi|\Psi_{01}> = |\Psi_{01}>$ then $|t> \in |\Psi_{00}>$, and $|t_2> = |0>$.

Again, it is important to note that in order to create register, $|\Psi_{00}>$, one first creates the state, $|0>^{(\otimes(n))}$ and then operates on this state the operator, $I \otimes I \otimes H^{(\otimes(n-2))}$ where *I* denotes the identity operator and *H* denotes the Hadamard operator. Therefore, if one will apply once more the same operator, $I \otimes I \otimes H^{(\otimes(n-2))}$, on thus formed register, $|0>^{(\otimes(n))}$ then one will get back the original state, $|0>^{(\otimes(n))}$ where one has started. Now, when Case(i) is true, i.e. when j$|\Psi_{00}>$ remains invariant under the action of the operator, $U_\Phi$ i.e. when $U_\Phi|\Psi_{00}> = |\Psi_{00}>$ then it is clear to see that one will get back the original starting state, , $|0>^{(\otimes(n))}$ when one will operate the same above mentioned operator, $I \otimes I \otimes H^{(\otimes(n-2))}$, on $U_\Phi|\Psi_{00}>$ and this will verify that Case (i) is true.

5. Suppose case (i) given above in 3 is true. Then $t \in [2(n-1), (2n) - 1]$, we then divide the register $|\Psi_1>$ into two equal parts, i.e. we create two new registers as follows: We construct two new quantum registers $|\Psi_{10}>$ and $|\Psi_{11}>$ such that $|\Psi_{10}>$ will represent the quantum register that contains computational basis states corresponding to all those indices, x, such that $2(n-1) \leq x \leq 2(n-1)+(2(n-2)-1)$ and $|\Psi_{11}>$ will represent the quantum register that contains computational basis states corresponding to all those indices, x, such that $(2(n-1) + 2(n-2)) \leq x \leq (2n - 1)$. Thus, we have
$|\Psi_{10}> = |1> \otimes |0> \otimes H^{(\otimes(n-2))}|0>^{(\otimes(n-2))}$
and,
$|\Psi_{11}> = |1> \otimes |1> \otimes H^{(\otimes(n-2))}|0>^{(\otimes(n-2))}$
Now, clearly: either $|t> \in |\Psi_{10}>$ or $|t> \in |\Psi_{11}>$, and this is to be determined in exactly similar way as described above 3 by doing the action of the the unitary operator, $U_\Phi$, on quantum registers $|\Psi_{10}>$ and $|\Psi_{11}>$, which will lead to creation of fully "distinguishable" states, $U_\Phi|\Psi_{10}>$ and $U_\Phi|\Psi_{11}>$.
Thus, we have two cases again:
Case (i): Without loss of generality, if $U_\Phi|\Psi_{10}> = |\Psi_{10}>$ then $|t> \in |\Psi_{11}>$, and $|t_2> = |1>$.
Case (ii): Without loss of generality, if $U_\Phi|\Psi_{11}> = |\Psi_{11}>$ then $|t> \in |\Psi_{10}>$, and $|t_2> = |0>$.

Again, it is important to note that in order to create register, $|\Psi_{10}>$, one first creates the state, $|0>^{(\otimes(n))}$, and then operates on this state the operator, $X\otimes I\otimes H^{(\otimes(n-2))}$ , where $I$ denotes the identity operator, $X$ denotes the NOT operator and $H$ denotes the Hadamard operator. Therefore, if one will apply once more the same operator, $X\otimes I\otimes H^{(\otimes(n-2))}$, on thus formed register, $|0>^{(\otimes(n))}$ then one will get back the original state, $|0>^{(\otimes(n))}$, where one has started. Now, when Case(i) is true, i.e. when $|\Psi_{10}>$ remains invariant under the action of the operator, $U_\Phi$ i.e. when $U_\Phi|\Psi_{10}> = |\Psi_{10}>$ then it is clear to see that one will get back the original starting state, $|0>^{(\otimes(n))}$, when one will operate the same above mentioned operator, $X\otimes I\otimes H^{(\otimes(n-2))}$, on $U_\Phi|\Psi_{10}>$ and this will verify that Case (i) is true.

On similar lines one can note that in order to create register, $|0>^{(\otimes(n))}$one first creates the state, $|0>^{(\otimes(n))}$, and then operates on this state the operator, $X\otimes X\otimes H^{(\otimes(n-2))}$ , where $X$ denotes the NOT operator and $H$ denotes the Hadamard operator. Therefore, if one will apply once more the same operator, $X\otimes X\otimes H^{(\otimes(n-2))}$ , on thus formed register, $|\Psi_{11}>$, then one will get back the original state, $|0>^{(\otimes(n))}$, where one has started. Now, when Case (ii) is true, i.e. when $|\Psi_{11}>$, remains invariant under the action of the operator, $U_\Phi$ i.e. when $U_\Phi|\Psi_{11}> = |\Psi_{11}>$ then it is clear to see that one will get back the original starting state, $|0>^{(\otimes(n))}$, when one will operate the same above mentioned operator, $X\otimes X\otimes H^{(\otimes(n-2))}$, on $U_\Phi|\Psi_{11}>$ and this will verify that Case (ii) is true.

6. We continue on these lines with dividing each time that (correct) quantum register to which the target state, |ti, belongs and this is achieved by creating the corresponding two appropriate new quantum registers as described above. It is clear to see that by carrying out these steps the size of the quantum register that contains the target state, |ti, will go on reducing and finally it will this register will be containing just one computational basis state, namely, the target state $|t>$ itself!

## 5.3. Concluding Remarks:

**Remark 5.3.1:** We divide during each iteration the quantum register to which the target state, $|t>$, belongs. Suppose at an intermediate stage we have arrived at the quantum register, say,

$$|\Psi i_1 i_2 \ldots\ i_k> = |i_1 i_2 \ldots\ i_k> \otimes H^{(\otimes(n-k))} |0>^{(\otimes(n-k))}$$

to which the target state, $|t>$, belongs. We then divide this quantum register into two new quantum registers by creating two new quantum registers as follows:

$$|\Psi i_1 i_2 \ldots\ i_k 0> = |i_1 i_2 \ldots\ i_k> \otimes |0> \otimes H^{(\otimes(n-k-1))} |0>^{(\otimes(n-k-1))}$$

and,

$$|\Psi i_1 i_2 \ldots\ i_k 1> = |i_1 i_2 \ldots\ i_k> \otimes |1> \otimes H^{(\otimes(n-k-1))} |0>^{(\otimes(n-k-1))}$$

Out of these above two newly created quantum registers again there will be only (some) one quantum register to which the target state, $|t>$, will belong, and so on. We continue on these lines till we finally reach the target state $|t>$ itself!

**Remark 5.3.2:** If we want to find the the correct answer to a K-SAT problem, containing 50 variables which has only one correct answer, by trial and error method then we may require to check $2^{50}$ = 1,125,699,906,842,624 possible answers, i.e. we will require to carry out in the worst case **$2^{50}$ = 1,125,699,906,842,624 answer checking operations.** On the other hand, using the above developed efficient quantum algorithm, we can directly reach to the correct answer in only 50 iterative operations, thus, in using this new proposed quantum algorithm we will require to carry out **only 50 answer building operations!**

## 5.4. New Efficient Quantum Algorithms for Solving NP-complete Problems:

As mentioned above in the section describing background of the invention, any NP-Complete problem can be mapped into any other NP-complete problem using only polynomial resources [12]. Thus, if an efficient algorithm were found that can solve one type of NP-Complete problem efficiently, this would immediately lead to efficient algorithms for all the other NP-complete problems (up to the polynomial cost of translation). We have developed in the above section an efficient quantum algorithm for solving a K-SAT, problem for $K \geq 3$, which is an NP-Complete problem, in fact the first problem shown to be NP-Complete [1], therefore, by mapping (one by one) every other NP-Complete problem to K-SAT problem with polynomial const of translation one can develop an efficient quantum algorithm for every other NP-Complete problem and thus can have an efficient quantum solution for each and every other NP-Complete problem.

## 5.5 An Illustrative Example:

To illustrate our approach we discuss steps to find the solution to a K–SAT problem, $K \geq 3$, containing four variables and defined by the Boolean function: $\Phi(x1, x2, x3, x4)$. Suppose there is a unique solution, t = (t1, t2, t3, t4) = (1, 1, 0, 0). We prepare two registers:

$|\Psi_0\rangle = |0\rangle \otimes H^{(3)}|000\rangle$

and

$|\Psi_1\rangle = |1\rangle \otimes H^{(3)}|000\rangle$

We will find that $U_\Phi|\Psi_0\rangle = |\Psi_0\rangle$. This implies that $|t\rangle \in |\Psi_1\rangle$ and since all the states present in the superposition of states belonging to this register $|\Psi_1\rangle$ begin with the single qubit state $|1\rangle$, therefore, $|t_1\rangle = |1\rangle$. Now, since $|t\rangle \in |\Psi_1\rangle$ therefore as per algorithm we now divide the register $|\Psi_1\rangle$ into two new registers. In other words, we now create two new registers (which is equivalent to breaking the register $|\Psi_1\rangle$ into two equal parts) as follows:

$|\Psi_{10}> = |1> \otimes |0> \otimes H^{(2)}|00>$

and

$|\Psi_{11}> = |1> \otimes |1> \otimes H^{(2)}|00>$

We will now find that $U_\Phi|\Psi_{10}> = |\Psi_{10}>$, therefore, $|t> \in |\Psi_{11}>$. Since the second qubit of all the states present in the superposition of states belonging to the register $|\Psi_{11}>$ is equal to $|1>$ therefore we have $|t_2> = |1>$. As per algorithm we now create two new registers (which is equivalent to breaking the register $|\Psi_{11}>$ into two equal parts) as follows:

$|\Psi_{110}> = |1> \otimes |1> \otimes |0> \otimes H^{(1)}|0>$

and

$|\Psi_{111}> = |1> \otimes |1> \otimes |1> \otimes H^{(1)}|0>$

We will now find that $U_\Phi|\Psi_{111}> = |\Psi_{111}>$ therefore $|t> \in |\Psi_{110}>$. Since the third qubit of all the states present in the superposition of states belonging to the register $|\Psi_{110}>$ is equal to $|0>$ therefore we have $|t_3> = |0>$. As per algorithm we now create two new registers (which is equivalent to breaking the register $|\Psi_{110}>$ into two equal parts) as follows:

$|\Psi_{1100}> = |1> \otimes |1> \otimes |0> \otimes |0>$

and,

$|\Psi_{1101}> = |1> \otimes |1> \otimes |0> \otimes |1>$

We will now find that $U_\Phi|\Psi_{1101}> = |\Psi_{1101}>$ therefore $|t> \in |\Psi_{1100}>$.

Now, there is only one state present in register $|\Psi_{1100}>$ and its fourth qubit is equal to $|0>$ therefore we have $|t_4> = |0>$.

Thus we have obtained

$|t> = |1> \otimes |1> \otimes |0> \otimes |0>$.

Thus, we determine the computational basis state |t> = |1100>, corresponding to the assignment t = (1, 1, 0, 0), representing the only solution for the given K-SAT problem in 4 iterative steps, i.e. we determine the computational basis state |t> = |1100>, corresponding to the assignment t = (1, 1, 0, 0), representing the only solution for the given K-SAT problem in only 4 answer building operations!